\documentclass[a4paper,11pt]{article} 
\usepackage[catalan,spanish,english]{babel} 
\usepackage[psamsfonts]{amssymb} 
\usepackage{cmmib57} 
\usepackage{exscale} 
\usepackage{amsmath} 
\usepackage{amscd} 
\usepackage{latexsym} 
\usepackage{graphicx} 
\def\1{{\rm 1 }\hskip -0.21truecm 1}

\newcommand{\Gam}{\Gamma} 
\newcommand{\gam}{\gamma} 
\newcommand{\ffi}{\varphi}

\newcommand{\red}{\mathbb{R}^{d}}

\newcommand{\sig}{\sigma}

\newcommand{\tf}{\mathcal{F}} 
\newcommand{\beq}{\begin{equation}} 
\newcommand{\eeq}{\end{equation}} 
\newcommand{\beqn}{\begin{equation*}} 
\newcommand{\eeqn}{\end{equation*}}

\newcommand{\qed}{\quad{$\square$}} 
\newtheorem{prop}{{\bf{Proposition}} } 
\newtheorem{obs}{{\bf{Remark}} } 
 
\newtheorem{corol}{{\bf{Corollary}} } 
\newtheorem{teor}{{\bf{Theorem}} } 
 
\newtheorem{ex}{\bf{Example}} 
\newtheorem{remark}{\bf{Remark}} 

\begin{document} 

%%%% TITLE PAGE
\begin{titlepage}
\null
\vspace{2cm}
\begin{center}
{\LARGE \bf Regularity of the sample paths of a class\\[2mm]
of second order spde's}\\[2mm]

\bigskip

by\\
\vspace{7mm}

\begin{tabular}{l@{\hspace{10mm}}l@{\hspace{10mm}}l}
{\sc Robert C. Dalang}$\,^{(\ast)}$  &and&{\sc Marta Sanz-Sol\'e}$\,^{(\ast\ast)}$\\
{\small Institut de Math\'ematiques}     &&{\small Facultat de Matem\`atiques}\\
{\small Ecole Polytechnique F\'ed\'erale}          &&{\small Universitat de Barcelona}\\
{\small 1015 Lausanne }            &&{\small Gran Via 585}\\         
{\small Switzerland}                 &&{\small 08007 Barcelona, Spain}\\
{\small e-mail: robert.dalang@epfl.ch}      &&{\small e-mail: marta.sanz@ub.edu}\\
\null

\end{tabular}
\end{center}

\vspace{1 cm}

{\bf Abstract:} We study the sample path regularity of the solutions of a class of
spde's which are second order in time and that includes the stochastic wave equation. Non-integer powers of the spatial
Laplacian are allowed.
The driving noise is white in time and spatially homogeneous. Continuing with the work initiated in Dalang and Mueller (2003),
we prove that the solutions belong to a fractional $L^2$-Sobolev space. We also prove H\"older continuity in time and
therefore, we obtain joint H\"older continuity in the time and space variables.
Our conclusions rely on a precise analysis of the properties of the stochastic integral used in the rigourous formulation of
the spde, as introduced by Dalang and Mueller. For spatial covariances given by Riesz kernels, we show that our results are
optimal. 

\bigskip


\medskip


\vspace{0.5 cm}

\footnotesize
{\begin{itemize}
\item[$^{(\ast)}$] Partially supported by the Swiss National Foundation for Scientific Research.
\item[$^{(\ast\ast)}$] Partially supported by the grant BFM2003-01345 from \textit{Direcci\'on
General de Investigaci\'on, Ministerio de Ciencia y Tecnolog\'{\i}a}, Spain.
\end{itemize}}

\end{titlepage}

%%%%END TITLE PAGE
\section{\bf Introduction} 

Modelling with stochastic partial differential equations (spde's) provides successful
understanding of the evolution of many physical phenomena. 
 A basic issue 
is how to choose the ingredients so that the spde possesses a solution in 
a strong sense---giving rise to a function-valued stochastic process---and to 
fix, in the most precise way possible, the 
function space that contains the sample paths of the solution. It is well known that 
this amounts to finding the right balance between the roughness of the driving noise---the stochastic input in the
model---and the singularities of the differential operator that defines the equation, 
which may depend on the dimension. 

This paper focusses on the analysis of the following class of spde's: 
\begin{align} 
&\left(\frac{\partial^2}{\partial t^2} + (- \Delta)^{(k)}\right) u(t,x) = \sigma\big(u(t,x)\big) \dot F(t,x) 
+ b \big(u(t,x)\big),\nonumber\\ 
&u(0,x) = v_0(x),\qquad 
 \frac{\partial}{\partial t}u(0,x) = \tilde v_0(x). \label{0.1} 
\end{align} 
In this equation, $t\in[0,T]$ for some fixed $T>0$, $x\in\red$, $d\in\mathbb{N}$, $k\in\,]0,\infty[$ and 
$\Delta^{(k)}$ denotes the fractional Laplacian on $\red$. This includes for instance the stochastic wave equation 
in any spatial dimension $d$. 
The coefficients $\sigma$ and $b$ are Lipschitz continuous functions 
and satisfy $|\sigma(z)| + |b(z)| \le C |z|$, for some positive constant $C$. 
The generalized process $\dot F$ is a Gaussian random field, 
white in time and spatially homogeneous with spatial correlation. More precisely, let $\Gam$ be a non-negative and
non-negative 
definite tempered measure on $\red$. Let $\mathcal{D}(\mathbb{R}^{d+1})$ be the space of Schwartz test functions 
(see \cite{schwartz}). On a probability space $(\Omega,\tf,P)$, we define a Gaussian process 
$F = \big(F(\ffi),\ \ffi\in \mathcal{D}(\mathbb{R}^{d+1})\big)$ with mean zero and covariance functional given by 
\beqn 
E\big(F(\ffi) F(\psi)\big) = \int_{\mathbb{R}_+} ds \int_{\red} \Gam(dx) (\ffi(s)*\tilde\psi(s))(x), 
\eeqn 
where $\tilde\psi(s)(x)=\psi(s)(-x)$. 

   Using an extension of Walsh's stochastic integral with respect to martingale measures \cite{walsh}, developed in
\cite{dalangmueller}, we give a rigourous meaning to problem (\ref{0.1}) in a mild form (see Equation (\ref{3.3})). In fact,
in \cite{dalangmueller} a particular case of Equation (\ref{0.1}) (when $k\in\mathbb{N}$ is an integer
and $b\equiv 0$) was introduced and studied. 

Let $\mu = {\tf}^{-1}\Gam$ be the spectral measure of $F$ and assume that 
\beq 
\label{0.2} 
\int_{\red} \frac{\mu(d\xi)}{(1 + |\xi|^2)^k} < \infty. 
\eeq 

Under suitable assumptions on the initial 
condition and the restrictions on $k$ and $b$ mentioned above, Theorem 9 in \cite{dalangmueller} establishes the existence
of a unique solution satisfying 
\beqn 
\sup_{0\le t\le T}E\big(\Vert u(t)\Vert_{L^2(\red)}^2\big) < \infty, 
\eeqn 
for which $t \mapsto u(t)\in L^2(\red)$ is mean-square continuous. 

   Here, we want  to study the regularity properties of the 
sample paths---both in time and in space---of Equation (\ref{0.1}), when conditions stronger than condition (\ref{0.2}) are
imposed. 

   In several examples of spde's driven by spatially homogeneous noise, one can prove that their solutions are 
real-valued random fields $u = \big(u(t,x)$, $(t,x)\in[0,T]\times \red\big)$. 
This is the case for the stochastic heat equation in any 
spatial dimension $d\ge 1$, for the stochastic wave equation in 
dimension $d\in\{1,2\}$, or even in dimension $d=3$, if the initial conditions 
vanish (\cite{dalangfrangos}, \cite{dalang}, \cite{milletss}). 
Joint H\"older continuity in $(t,x)$ of the sample paths of the solution can usually be obtained using Kolmogorov's 
continuity condition (\cite{milletss}, \cite{sssarra}). However, for more general equations, such as those considered in this
paper, 
one can only expect solutions $u = \big(u(t),\ t\in[0,T]\big)$ taking values 
in some function space. Kolmogorov's condition is still well suited for establishing 
regularity properties in time, but it is not for the study of spatial continuity. Regularity in space may be obtained by means
of
Sobolev  type imbeddings, if one could prove that the solution takes values in some fractional Sobolev 
space $H_p^{\alpha}$, $p\in[1,\infty[$, $\alpha\in [0,\infty[$. Indeed, $H_p^{\alpha}$ is 
imbedded in the space of $\gamma$-H\"older continuous functions $\mathcal{C}^\gam(\red)$, for 
any $\gam\in\,]0,\alpha-\frac{d}{p}[$, whenever $\alpha > \frac{d}{p}$. This fact explains one of the main
advantages of an 
$L^p$-theory for spde's, for arbitrary values of $p$, leading to optimal results in $\gamma$. 

   Until now, $L^p$-theory for spde's has been mainly 
developed for parabolic spde's (see for instance \cite{krylov} and the references herein). 
Recently, we have been able to use an $L^p$ approach to study the sample path behaviour in $(t,x)$ of the stochastic wave
equation in dimension $d=3$ (see \cite{dalangss}), driven by the type of noise described above and with a covariance
function whose singularity is given by a Riesz kernel. The methods used in the analysis of this particular equation are very
much related to the special form of the fundamental solution of the equation and of the covariance function of the noise;
they do not seem to be exportable to the more general situation we are considering here. 

In this paper, we establish sufficient conditions on the spectral measure $\mu$ of the noise that ensure that the solution
of Equation (\ref{0.1}) belongs a.s.~to some fractional Sobolev space $H_2^{\alpha}$, for some $\alpha\in [0,k[$. Then we
prove H\"older continuity in time of the solution and show that the results are optimal when the covariance measure is a
Riesz kernel. 

Let $M$ be the martingale measure extension of the process $F$, obtained in \cite{dalangfrangos} 
(see also \cite{dalang}) and let $Z$ be an $L^2(\red)$-valued stochastic process. In Section \ref{sec2}, we prove that,
under suitable assumptions on the 
$\mathcal{S}'(\red)$-valued function $G$, the stochastic integral 
\beqn 
v_{G,Z}(T) = \int_0^T \int_{\red}G(s,\cdot-y)Z(s,y) M(ds,dy) 
\eeqn 
introduced in \cite{dalangmueller}, which defines a random element of $L^2(\red)$, belongs
 in fact to $H_2^{\alpha}(\red)$ and is such that $E\big(\Vert v_{G,Z}(T)\Vert^2_{H_2^{\alpha}(\red)}\big)
< \infty$. To establish this fact, we will need to prove the existence of the Fourier
transform of the stochastic integral 
$v_{G,Z}(T)$.  Recall \cite{treves} that for a function $g \in H_2^{\alpha}(\red)$, 
\beqn
   \Vert g\Vert^2_{H_2^{\alpha}(\red)} = \int_{\red} d\xi \, (1 + |\xi|^2)^\alpha \, |\tf g(\xi)|^2,
\eeqn
where, for $\varphi \in C^\infty_0(\red)$,
\beqn
   \tf \varphi(\xi) = \int_{\red} dx\, e^{i\xi x} \, \varphi(x).
\eeqn

   Let $\mathcal{L} = \partial_{tt}^2 + (- \Delta)^{(k)}$, $k\in\,]0,\infty[$. We prove that, if for some $\alpha\in [0,k[$,
\beq 
\label{0.3} 
   \int_{\red} \frac{\mu(d\xi)}{(1+|\xi|^2)^{k-\alpha}} < \infty,
\eeq 
then the preceding result applies to the fundamental solution of $\mathcal{L} u = 0$. 
For the other results of this paper, we will assume property (\ref{0.3}). We note that we treat indifferently the case of
integer and fractional powers of the Laplacian.

   Section \ref{sec2} is devoted to studying path properties in time of the stochastic integral 
\beqn 
v_{G,Z}(t) = \int_0^t \int_{\red}G(s,\cdot-y)Z(s,y) M(ds,dy)
\eeqn 
and the H\"older continuity of
\beq 
   u_{G,Z}(t) = \int_0^t \int_{\red}G(t-s,\cdot-y)Z(s,y) M(ds,dy).
\label{0.3.1}
\eeq
We first identify the increasing process of the $H_2^{\alpha}(\red)$-valued martingale 
$\big(v_{G,Z}(t), t\in[0,T]\big)$. Fix $\alpha\in [0,k[$ and assume that there exists 
$\eta\in\,]\frac{\alpha}{d},1[$ such that the following condition, which is stronger than (\ref{0.3}), holds: 
\beq 
\label{0.4} 
\int_{\red} \frac{\mu(d\xi)}{(1+|\xi|^2)^{k\eta-\alpha}} < \infty. 
\eeq 
Using Kolmogorov's continuity condition, we obtain that the sample paths of $(u_{G,Z}(t),\  t\in[0,T])$ are a.s.~H\"older
continuous. In the particular case where $\Gam (dx) = |x|^{-\beta}$, $\beta\in\,]0,d[$, the results are proved to be
optimal. By means of the Sobolev imbedding theorem, we also obtain  H\"older continuity in the space variable. However,
the conditions for validity of this result are rather restrictive. 

In Section 3, we transfer the results of the preceding sections to the solution of Equation (\ref{0.1}). 
Fix $\alpha\in [0,k[$, assume (\ref{0.3}) and that the initial conditions are such that $v_0\in H_2^{\alpha}(\red)$ and 
$\tilde v_0\in H_2^{\alpha - k}(\red)$. We prove the existence of a solution to (\ref{0.1}) satisfying 
\beqn 
\sup_{0\le t\le T} E\big(\Vert u(t)\Vert_{H_2^{\alpha}(\red)}^q\big) < \infty. 
\eeqn 
Replacing assumption (\ref{0.3}) by (\ref{0.4}) and under additional (but natural) hypotheses on the initial conditions, we
obtain H\"older continuity in time of the solution of the equation. 

%%%%%%%%%%% END OF THE INTRODUCTION 

\section{\bf The stochastic integral as a random vector with values in a fractional Sobolev space}\label{sec2} 

   In this section, we consider the stochastic integral defined in Theorem 6 of \cite{dalangmueller}. Our aim 
is to prove that under suitable assumptions, this integral takes its values in the fractional Sobolev 
space $H_2^{\alpha}(\red)$, for some $\alpha\in\,]0,\infty[$. 

   Throughout this section, let ${\tf}_s$ be the $\sigma$-field generated by the martingale measure $(M_t,\ 0\le t\le s)$ 
described in the introduction. We consider a stochastic 
process $Z = (Z(s),\ s\in[0,T])$ with values in $L^2(\red)$ such that $Z(s)$ is ${\tf}_s$-measurable 
and the mapping $s\mapsto Z(s)$ is mean-square continuous from $[0,T]$ into $L^2(\red)$.

  The main result of this section is as follows. 

\begin{teor} 
\label{t1} 
Consider a deterministic map $G: [0,T] \longrightarrow {\mathcal S}'(\red)$. Fix $\alpha\in [0,\infty[$
and assume that the following three conditions hold: 
\begin{description} 
\item{{\rm(i)}} For each $s\in[0,T]$, $\tf G(s)$ is a function and 
\beqn 
\sup_{0\le s\le T} \sup_{\xi\in\red} (1 + |\xi|^2)^{\frac{\alpha}{2}} \big|\tf G(s)(\xi)\big| < \infty.
\eeqn 
%\item{(ii)} For all $\psi\in \mathcal{S}(\red)$, 
\item{\rm{(ii)}} For all $\psi\in \mathcal{C}^{\infty}_0(\red)$, 
\beqn 
\sup_{0\le s\le T} \sup_{x\in\red}\big|\big(G(s) *\psi\big)(x)\big| < \infty.
\eeqn 
\item{{\rm(iii)}} 
\beqn 
\int_0^T ds \sup_{\xi\in\red} \int_{\red} \mu(d\eta) (1 + |\xi-\eta|^2)^\alpha \big|\tf G(s)(\xi - \eta)\big|^2 
<\infty. 
\eeqn 
\end{description} 
Then the stochastic integral 
\beqn
   v_{G,Z}(T) = \int_0^T \int_{\red} G(s,\cdot - y) Z(s,y) M(ds,dy)
\eeqn
satisfies 
\beqn 
E\big(\Vert v_{G,Z}(T)\Vert^2_{H_2^{\alpha}(\red)}\big) < \infty
\eeqn 
and 
\begin{align} 
E\big(\Vert v_{G,Z}(T)\Vert^2_{H_2^{\alpha}(\red)}\big) 
&= E\big(\Vert v_{(I-\Delta)^{\frac{\alpha}{2}}G,Z}(T)\Vert^2_{L^2(\red)}\big) \nonumber\\ 
&= I^{\alpha}_{G,Z}, \label{1.0} 
\end{align} 
where 
\beqn 
v_{(I-\Delta)^{\frac{\alpha}{2}}G,Z}(T) = \int_0^T \int_{\red} (I-\Delta)^{\frac{\alpha}{2}}G(s,\cdot-y) 
Z(s,y) M(ds,dy)
\eeqn
and
\begin{align}
I^{\alpha}_{G,Z} &= \int_0^Tds \int_{\red} d\xi E\big(|\tf Z(s)(\xi)|^2\big)\nonumber \\ 
& \qquad \times\int_{\red}\mu(d\eta) (1 + |\xi-\eta|^2)^{\alpha} \big|\tf G(s)(\xi - \eta)\big|^2. \label{1.1} 
\end{align} 
\end{teor} 

The proof of this theorem relies on a preliminary result that identifies the Fourier transform of the stochastic integral 
$v_{G,Z}$ for $G$ and $Z$ satisfying more restrictive assumptions than those above, namely:
\medskip

\begin{description}
\item{\bf (G1')} For each $s\in[0,T]$, $G(s) \in \mathcal{C}^{\infty}(\red)$, $\tf G(s)$ is a function, 
\beqn 
\sup_{0\le s\le T}\sup_{x\in\red}|G(s,x)| < \infty
 \qquad\mbox{and}\qquad
\sup_{0\le s\le T}\sup_{\xi\in\red}|\tf G(s)(\xi)| < \infty. 
\eeqn 

\item{\bf(G2)} For $s\in[0,T]$, $Z(s) \in \mathcal{C}_0^\infty(\red)$ a.s., and there is a compact set 
$K\subset \red$ such that ${\rm supp}\, Z(s) \subset K$, for $s\in[0,T]$. In addition, 
the mapping $s\mapsto Z(s)$ is mean-square continuous from $[0,T]$ into $L^2(\red)$. 
\smallskip 

\item{\bf(G3)} $I_{G,Z} < \infty$, where 
\beq 
\label{1.2} 
I_{G,Z} = \int_0^T ds \int_{\red} d\xi\, E(|\tf Z(s)(\xi)|^2) \int_{\red} \mu(d\eta) |\tf G(s)(\xi-\eta)|^2. 
\eeq 
\end{description}
Notice that our assumption {\bf (G1')} is stronger than {\bf (G1)} in \cite{dalangmueller} (which does not 
suppose the boundedness of the Fourier transform of $G$), while 
{\bf(G2)} and {\bf(G3)} appear with the same name in \cite{dalangmueller}. 

Under {\bf (G1')}, {\bf(G2)}, and {\bf(G3)}, the stochastic integral 
\beqn 
v_{G,Z}(T)(x) = \int_0^T \int_{\red} G(s,x-y) Z(s,y) M(ds,dy) 
\eeqn 
is well-defined, for any $x\in\red$, as a Walsh stochastic integral (see Lemma 1 in \cite{dalangmueller}). 
The integral 
\beqn 
\int_0^T \int_{\red}\tf G(s, \cdot - y)(\xi) Z(s,y) M(ds, dy) 
\eeqn 
is also well-defined as a Walsh stochastic integral.
 Indeed, $\tf G(s)(\cdot - y)(\xi) = e^{i\xi\cdot y} \tf G(s)(\xi)$, and
\begin{align*} 
&E\Big(\int_0^T ds \int_{\red} \Gam (dy) \int_{\red} dz |e^{i\xi\cdot z} \tf G(s)(\xi) Z(s,z)
%\\ &\qquad\qquad\qquad \times 
e^{i\xi (y-z)} \tf G(s)(\xi) Z(s,y-z)|\Big)\\ 
&\qquad\le \sup_{0\le s\le T}\sup_{\xi\in\red}|\tf G(s)(\xi)|^2 \int_0^T ds \int_{\red} \Gam (dy) 
E\big(|Z(s,\cdot)| \ast |\tilde Z(s,\cdot)|\big)(y)\\ 
&\qquad \le C \int_0^T ds\, E\big(\Vert Z(s,\cdot)\Vert^2_{L^2(\red)}\big) \Gam(K-K) < \infty. 
\end{align*}

\begin{prop} 
\label{p1} 
We assume the hypotheses ${\bf (G1')}$, ${\bf (G2)}$ and ${\bf (G3)}$.
Then the Fourier transform $\tf v_{G,Z}(T)$ of the stochastic integral $v_{G,Z}(T)$ is given by 
\beqn 
\tf v_{G,Z}(T)(\xi)=\int_0^T \int_{\red} \tf G(s,\cdot - y)(\xi) Z(s,y) M(ds, dy). 
\eeqn 
\end{prop} 

\noindent{\em Proof}. Let $\ffi \in \mathcal{S}(\red)$. We want to check that 
\beq 
\label{1.3} 
\langle v_{G,Z}(T), {\tf}^{-1}\ffi \rangle = \langle \int_0^T \int_{\red} \tf G(s,\cdot - y)(\cdot) Z(s,y) M(ds, dy), 
\ffi\rangle, 
\eeq 
where $\langle \cdot,\cdot\rangle$ denotes the inner product in $L^2(\red)$. 

We verify the assumptions of the stochastic Fubini's theorem in \cite{walsh}: since $G$ is uniformly bounded and $Z(s)$
has compact support, 
\begin{align*} 
&E\Big ( \int_{\red} dx \int_0^T ds \int_{\red} \Gam(dy)\int_{\red} dz\, |{\tf}^{-1}\ffi (x)|^2 
|G(s,x-z)|\\ 
&\qquad\qquad \times |Z(s,z)|\, |G(s,x-z+y)|\, |Z(s,z-y)|\Big )\\ 
&\le C E\Big ( \int_{\red} dx\, |{\tf}^{-1}\ffi (x)|^2 \int_0^T ds \int_{\red}\Gam(dy)\int_{\red} dz \,
|Z(s,z)|\, |\tilde Z(s,y-z)|\Big )\\ 
&\leq C \Vert \ffi\Vert^2_{L^2(\red)} \Gam(K-K) \int_0^T ds\, E\big(\Vert Z(s)\Vert^2_{ L^2(\red)}\big) < \infty. 
\end{align*} 
Applying this Fubini's theorem and Plancherel's identity, we obtain 
\begin{align} \nonumber
\langle v_{G,Z}(T), {\tf}^{-1}\ffi \rangle &= \int_{\red} dx\, \big(\int_0^T \int_{\red} {\tf}^{-1} \ffi (x) 
G(s,x-y) Z(s,y) M(ds,dy)\big)\\ \nonumber
& = \int_0^T \int_{\red} \big(\int_{\red} dx\, {\tf}^{-1} \ffi (x) G(s,x-y)\big) Z(s,y) M(ds,dy)\big)\\ 
& = \int_0^T \int_{\red} \big(\int_{\red} d\xi\, \ffi(\xi) \tf G(s)(\xi) \exp(i y\cdot\xi)\big) 
Z(s,y) M(ds,dy). 
\end{align} 
In order to apply again the stochastic Fubini's theorem, we note that 
\begin{align} \nonumber
&E\Big( \int_{\red} d \xi \int_0^T ds \int_{\red} \Gam(dy)\int_{\red} dz\, |\ffi(\xi)|^2 |\tf G(s)(\xi)|^2 
|Z(s,z)|\, |Z(s,y-z)|\Big)\\ \nonumber
&\le\Gam(K-K) \int_{\red} d\xi \int_0^T ds\, |\ffi(\xi)|^2 |\tf G(s)(\xi)|^2 E\big(\Vert Z(s)\Vert^2_{ L^2(\red)}\big)\\ 
&\le C \Vert\ffi\Vert^2_{ L^2(\red)}\sup_{0\le s\le T}E\big(\Vert Z(s)\Vert^2_{ L^2(\red)}\big) 
\sup_{0\le s\le T}\sup_{\xi\in\red}|\tf G(s)(\xi)|^2. 
\label{rde8a}
\end{align} 
Therefore, applying again the above-mentioned Fubini's theorem shows that the last right-hand side of (\ref{rde8a}) is
equal to 
\beqn 
%\langle v_{G,Z}(T), {\tf}^{-1}\ffi \rangle =
   \int_{\red} d\xi\, \ffi(\xi) \Big( \int_0^T \int_{\red} \tf G(s,\cdot - y)(\xi) Z(s,y) M(ds,dy)\Big),
\eeqn 
which establishes (\ref{1.3}). \hfill \qed 

\bigskip 

\noindent{\em Proof of Theorem \ref{t1}}. We proceed in several steps. 

{\it Step 1.} Assume first that $G$ and $Z$ satisfy the assumptions {\bf (G1')}, {\bf (G2)} and {\bf (G3)}. 
Suppose also that $G^{\alpha}(s):= (I-\Delta)^{\frac{\alpha}{2}}G(s)$ satisfies 
{\bf (G1')} and that $I_{G,Z}^{\alpha} < \infty$ (this last condition is implied by {\rm (iii)}). 
Then the stochastic integral 
$v_{(I-\Delta)^{\frac{\alpha}{2}}G,Z}(T)$ is well-defined in Walsh's sense and satisfies 
\beq 
\label{1.3.1} 
E\big(\Vert v_{(I-\Delta)^{\frac{\alpha}{2}}G,Z}(T)\Vert^2_{L^2(\red)}\big) = I_{G,Z}^{\alpha}. 
\eeq 
Indeed, this follows from Lemma 1 in \cite{dalangmueller}. 
Moreover, Proposition \ref{p1} implies that
\beq\label{rde8b}
\tf v_{(I-\Delta)^{\frac{\alpha}{2}}G,Z}(T) = v_{\tf (I-\Delta)^{\frac{\alpha}{2}} G,Z}(T). 
\eeq

   By the definition of the norm in $H^{\alpha}_2(\red)$, Plancherel's theorem 
and Proposition \ref{p1}, we obtain 
\begin{align} 
& E\big( \Vert v_{G,Z}(T)\Vert^2_{H_2^{\alpha}(\red)}\big) = E\big(\int_{\red} d\xi\, (1+|\xi|^2)^{\alpha} 
\big| \tf v_{G,Z}(T)(\xi)\big|^2\big)\nonumber\\ 
& = E\big(\int_{\red} d\xi\, (1+|\xi|^2)^{\alpha} \big|v_{\tf G,Z}(T)(\xi)\big |^2\big) \nonumber\\ 
& = E\Big(\int_{\red} d\xi\, \big |\int_0^T\int_{\red} (1+|\xi|^2)^{\frac{\alpha}{2}} \tf G(s,\cdot - y)(\xi) 
Z(s,y) M(ds,dy)\big |^2 \Big) \nonumber\\ 
& = E\Big(\int_{\red} d\xi\, \big |\int_0^T\int_{\red} \tf \big( (I-\Delta)^{\frac{\alpha}{2}}G(s,\cdot - y))(\xi) 
Z(s,y) M(ds,dy)\big |^2 \Big)\nonumber\\ 
& = E\big(\Vert v_{(I-\Delta)^{\frac{\alpha}{2}}G,Z}(T)\Vert^2_{L^2(\red)}\big). \label{1.4} 
\end{align} 

Consequently, the theorem is proved in this particular situation. Notice that 
(\ref{1.0}) follows from (\ref{1.3.1}) and (\ref{1.4}). 
\medskip

{\it Step 2.} Assume  that $G(s)$ and $G^{\alpha}(s)$ satisfy {\bf (G1')} and that condition (iii) holds. By Lemma 3 in
\cite{dalangmueller}, there exists a sequence of stochastic processes 
$(Z_n, n\ge 1)$ satisfying 
{\bf (G2)} and {\bf (G3)} such that $\lim_{n\to\infty} I^0_{G^\alpha,Z_n-Z}=0$ and 
\beqn 
  v_{G^\alpha,Z}(T) = \lim_{n\to\infty} v_{G^\alpha,Z_n}(T), 
\eeqn 
with the limit taken in $L^2(\Omega; L^2(\red))$. The properties of $G^{\alpha}$ ensure that $\lim_{n\to\infty}
I^\alpha_{G,Z-Z_n} = 0$ as well. 

   We want to prove that $(v_{G,Z_n}(T),\ n\ge 1)$ is a Cauchy sequence in \break $L^2(\Omega;
H_2^{\alpha}(\red))$. Since 
$ H_2^{\alpha}(\red))$ is imbedded in $L^2(\red)$, the two limits of the sequence---in $L^2(\Omega; 
H_2^{\alpha}(\red))$ and in $L^2(\Omega;L^2(\red))$---must coincide. 

By the results proved in Step 1, 
\beq 
\lim_{n,m\to\infty}E\big(\Vert  v_{G,Z_n-Z_m}(T)\Vert^2_{H_2^{\alpha}(\red)}\big) 
= \lim_{n,m\to\infty} I^{\alpha}_{G,Z_n-Z_m} = 0. 
\eeq 
Let us now prove (\ref{1.0}) in this particular case. The previous convergence, the results 
stated in the first part of the proof and Lemma 3 in \cite{dalangmueller} applied 
to $g:= (I - \Delta)^{\frac{\alpha}{2}}G$ yield 
\begin{align*} 
E\big(\Vert v_{G,Z}(T)\Vert^2_{H_2^{\alpha}(\red)}\big)& = 
\lim_{n\to\infty} E\big(\Vert v_{G,Z_n}(T)\Vert^2_{H_2^{\alpha}(\red)}\big)\\ 
& = \lim_{n\to\infty} E\big(\Vert v_{(I - \Delta)^{\frac{\alpha}{2}}G,Z_n}(T)\Vert^2_{L^2(\red)}\big)\\ 
& = \lim_{n\to\infty} I^{\alpha}_{G,Z_n}\\ 
& = I^{\alpha}_{G,Z}. 
\end{align*} 
\medskip

{\it Step 3.} Let us now put ourselves under the assumptions of the theorem. 
Let $(\psi_n, n\ge 1)$ be an approximation of the identity such that $|\tf \psi_n(\xi)| \leq 1$, for all $\xi \in \red$. Set 
$G_n(s) = G(s) * \psi_n$, $G_n^{\alpha} = (I-\Delta)^{\frac{\alpha}{2}}G_n$. We now check that 
$G_n$ and $G_n^{\alpha}$ satisfy {\bf (G1')} and the assumption (iii) of the theorem. 

It is clear that $G_n(s)\in \mathcal{C}^{\infty}(\red)$. Moreover, condition (ii) yields 
\beqn 
\sup_{0\le s \le T} \sup_{x\in\red} |G_n(s,x)|) < \infty. 
\eeqn 
Since $\tf G(s)$ is a function, $\tf G_n(s)$ is also a function and (i) implies 
\beqn 
\sup_{0\le s \le T} \sup_{\xi\in\red} |\tf G_n(s)(\xi)| \le \sup_{0\le s \le T} \sup_{\xi\in\red}|\tf G(s)(\xi)| < \infty. 
\eeqn 
Notice that $G_n^{\alpha}(s)= G(s) * (I - \Delta)^{\frac{\alpha}{2}} \psi_n$. Since $(I -
\Delta)^{\frac{\alpha}{2}} 
\psi_n \in \mathcal{S}(\red)$, the Schwartz space of $C^\infty$ test functions with rapid decrease, we have
$G_n^{\alpha}\in \mathcal{C}^{\infty}(\red)$ (see for instance \cite{gasquetw}, Proposition 32.1.1). 

The condition $\sup_{0\le s \le T} \sup_{x\in\red}|G_n^{\alpha}(s,x)| < \infty$ is a consequence of assumption 
(ii). Since $\tf G(s)$ is a function, so is $\tf G_n^{\alpha}(s)$. Moreover, condition (i) yields 
\begin{align*} 
\sup_{n\ge 1}\sup_{0\le s \le T} \sup_{\xi\in\red} |\tf G^\alpha_n (s)(\xi)| 
& = \sup_{n\ge 1}\sup_{0\le s \le T} \sup_{\xi\in\red} |\tf G(s)(\xi)|\, |(1+|\xi|^2)^{\frac{\alpha}{2}} \tf \psi_n(\xi)|\\ 
& \le \sup_{0\le s \le T} \sup_{\xi\in\red}|(1+|\xi|^2)^{\frac{\alpha}{2}} |\tf G(s)(\xi)|\\
& < \infty. 
\end{align*} 

Consider the sequence of stochastic integrals $(v_{G_n,Z}(T), n\ge 1)$. Theorem 6 in \cite{dalangmueller} 
shows that $v_{G,Z}(T)$ is well-defined as an $L^2(\Omega;L^2(\red))$-valued random variable and 
\begin{align*} 
&E\big(\Vert v_{G_n,Z}(T) - v_{G,Z}(T)\Vert^2_{L^2(\red)}\big) = I_{G_n-G,Z}\\ 
& = \int_0^T ds \int_{\red} d \xi\, E(|\tf Z(s)(\xi)|^2) \int_{\red} \mu (d\eta) |\tf G(s)(\xi - \eta)|^2 |\psi_n(\xi) - 1|^2. 
\end{align*} 
By dominated convergence, this expression tends to zero as $n$ tends to infinity. 

We want to prove that $(v_{G_n,Z}(T), n\ge 1)$ is a Cauchy sequence in
\break
$L^2(\Omega;H_2^{\alpha}(\red))$.  Indeed, by the results stated in Step 2, we obtain 
\begin{align*} 
&E\big(\Vert v_{G_n-G_m,Z}(T)\Vert^2_{H_2^{\alpha}(\red)}\big) = I^{\alpha}_{G_n-G_m,Z}\\ 
& = \int_0^T ds \int_{\red} d \xi\, E(|\tf Z(s)(\xi))|^2) \\ 
&\quad \times \int_{\red} \mu(d\eta) (1+|\xi - \eta|^2)^{\alpha} 
|\tf G(s)(\xi)|\, |(\tf \psi_n - \tf \psi_m)(\xi)|^2. 
\end{align*} 
This last expression tends to zero as $n,m$ tend to infinity, by dominated convergence and 
assumption (iii). We have therefore proved that 
\beqn 
\lim_{n\to\infty}E\big(\Vert v_{G_n,Z}(T) - v_{G,Z}(T)\Vert^2_{H_2^{\alpha}(\red)}\big) = 0. 
\eeqn 

   By the results of Step 2, we obtain 
\beqn 
E\big(\Vert v_{G,Z}(T)\Vert^2_{H_2^{\alpha}(\red)}\big) 
= \lim_{n\to\infty}E\big(\Vert v_{G_n,Z}(T)\Vert^2_{H_2^{\alpha}(\red)}\big)\\ 
%& = \lim_{n\to\infty}E\big(||v_{(I - \Delta)^{\frac{\alpha}{2}}G_n,Z}^T||^2_{L^2(\red)}\big)\\ 
 = I^{\alpha}_{G,Z}. 
\eeqn 
This finishes the proof of the theorem.\hfill \qed 

\begin{ex} 
\label{e1} 

Consider the differential operator $\mathcal{L} = \partial_{tt}^2 + (- \Delta)^{(k)}$, $k\in\,]0,\infty[$, and denote by $G$
the fundamental solution of $\mathcal{L}u = 0$. It is easy to check that 
\beq\label{rde11a}
\tf G(t)(\xi) = \frac{\sin(t|\xi|^k)}{|\xi|^k}. 
\eeq
Fix $\alpha \in [0, k[$ and assume that 
\beq 
\label{1.5} 
\int_{\red} \frac{\mu (d \xi)}{(1 + |\xi|^2)^{k-\alpha}} < \infty. 
\eeq 
Then the assumptions of Theorem \ref{t1} are satisfied. 
\end{ex} 

Indeed, $\tf G(s)\in \mathcal{C}^{\infty}(\red)$. Moreover, 
\begin{align*} 
&\sup_{0\le s\le T} \sup_{\xi\in\red}\Big((1+|\xi|^2)^{\frac{\alpha}{2}} \frac{\sin(t|\xi|^k)}{|\xi|^k}\Big)\\ 
&\qquad\qquad \le T \sup_{|\xi|\le 1}(1+|\xi|^2)^{\frac{\alpha}{2}} + \sup_{|\xi|\ge 1}(1+|\xi|^2)^{\frac{\alpha-k}{2}}. 
\end{align*} 
Since $\alpha < k$, this last expression is finite and thus condition (i) is satisfied. 

   Let $\psi\in\mathcal{S}(\red)$. Then, 
\begin{align*} 
\sup_{0\le s\le T} \sup_{x\in\red} | \big(G(s)*\psi\big)(x) |& \le \sup_{0\le s\le T}\Vert\tf (G(s) *
\psi)\Vert_{L^1(\red)}\\ 
& \le T \Vert\tf \psi\Vert_{L^1(\red)} < \infty, 
\end{align*} 
proving (ii). 

Finally, we prove (iii). For any $k\in\,]0,\infty[$, it is easy to check that 
\beq\label{rde12a}
\sup_{0\le s\le T}|\tf G(s)(\xi)|^2 \le \frac{2^k(1+T^2)}{(1+|\xi|^2)^k}. 
\eeq
Therefore, 
\beq 
\label{1.5.1} 
\sup_{0\le s\le T}\int_{\red} \mu(d\eta) (1+|\xi - \eta|^2)^{\alpha} |\tf G(s)(\xi - \eta)|^2 \le C \int_{\red} 
\frac{\mu(d\eta)}{(1+|\xi - 
\eta|^2)^{k-\alpha}}, 
\eeq 
where $C$ is a positive constant depending on $T$ and $k$. 

Set $\gam=k-\alpha$. We will show that 
\beq 
\label{1.6} 
\sup_{\xi\in\red} \int_{\red} \frac{\mu(d\eta)}{(1+|\xi - \eta|^2)^{\gam}} 
\le \int_{\red} \frac{\mu(d\eta)}{(1+|\eta|^2)^{\gam}}. 
\eeq 
Combining this property with assumption (\ref{1.5}) yields (iii). Note for future reference that
\beq 
\label{1.6.1} 
\sup_{0\le s\le T}\sup_{\xi\in\red}\int_{\red} \mu (d\eta)(1+|\xi - \eta|^2)^{\alpha} |\tf G(s)(\xi - \eta)|^2 < \infty.
\eeq 

   In order to prove (\ref{1.6}), set $\tau_y(x)=x+y$. Following an argument that appears in \cite{KZ}, observe that
\begin{align*} 
&\int_{\red} \mu(d\eta)\, \frac{e^{-2\pi^2t|\eta|^2}}{(1+|\xi - \eta|^2)^{\gam}} = 
\int_{\red} \Gam (dx) {\tf}^{-1}\big( e^{-2\pi^2t|\cdot|^2}\tau_{-\xi}(1+|\cdot|^2)^{-\gam}\big)(x)\\ 
& = \int_{\red} \Gam (dx) \big({\tf}^{-1}( e^{-2\pi^2t|\cdot|^2})*{\tf}^{-1}(\tau_{-\xi}(1+|\cdot|^2)^{-\gam}) \big)(x)\\ 
& = \int_{\red} \Gam (dx)(p_t * e_{\xi} G_{d,\gam}(x), 
\end{align*} 
where $p_t = {\tf}^{-1}( e^{-2\pi^2t|\cdot|^2})$ is the Gaussian density with mean $0$ and variance $t$, $e_{\xi}(x) =
e^{2\pi i \langle x,\xi\rangle}$ and 
$ G_{d,\gam}(x) = {\tf}^{-1}(1+|\cdot|^2)^{-\gam}(x)$. 

Since both $p_t $ and $G_{d,\gam}$ are positive functions, 
\beqn 
   |(p_t * e_{\xi} G_{d,\gam}(x)| \le \int_{\red} p_t(y) G_{d,\gam}(x-y) dy. 
\eeqn 
Using monotone convergence and the fact that $p_t$ is a probability density function on $\red$, we obtain 
\begin{align*} 
\sup_{\xi\in \red}\int_{\red} \frac{\mu(d\eta)}{(1+|\xi - \eta|^2)^{\gam}} &= 
\sup_{\xi\in \red}\lim_{t\to 0}\int_{\red}\mu(d\eta) \frac{ e^{-2\pi^2t|\eta|^2}}{(1+|\xi - \eta|^2)^{\gam}}\\ 
& \le \lim_{t\to 0}\int_{\red} \Gam (dx)\int_{\red} G_{d,\gam}(x-y) p_t(y) dy\\ 
& \le \sup_{y\in\red} \int_{\red} \Gam (dx) G_{d,\gam}(x-y). 
\end{align*} 
However, 
\beq 
\sup_{y\in\red} \int_{\red} \Gam (dx) G_{d,\gam}(x-y) = \int_{\red}\frac{\mu(d\xi)}{(1+|\xi|^2)^{\gam}}, 
\eeq 
(see for instance \cite{sanzsole}, Lemma 9.8). This proves (\ref{1.6}). 
\medskip 

\begin{ex} 
\label{e2} 
Let $\mathcal{L}$ and $G$ be as in Example \ref{e1}. Assume that $\Gam(dx) = |x|^{-\beta}$, with 
$\beta\in \,]0,d[$. Then $\mu(d\xi) = C |\xi|^{-d+\beta}$ (see for instance \cite{donoghue}). Elementary 
calculations show that (\ref{1.5}) holds provided that $\beta < 2k$ and $\alpha\in [0,k-\frac{\beta}{2}[$. 
\end{ex} 
\medskip 

Fix $q\in\,]1,\infty[$ and $s\in\,]0,\infty[$. It is well known that the fractional Sobolev space $H_q^s(\red)$, 
is imbedded in the space $\mathcal{C}^{\gam}(\red)$ of $\gam$-H\"{o}lder continuous functions with $\gam \le
s-\frac{d}{q}$, 
whenever $s-\frac{d}{q} > 0$. Moreover, if $1<q<d$, $d>sq$, then $H_q^s(\red)$ 
is imbedded in $L^p(\red)$, for $q < p < \frac{dq}{d-sq}$ 
(see \cite{shimakura}). This yields the following. 

\begin{corol} 
\label{c1} 
\begin{enumerate} 
\item 
Suppose that the assumptions of Theorem \ref{t1} are satisfied for some $\alpha \in\,]\frac{d}{2}, \infty[$. 
Then almost surely, the stochastic integral $v_{G,Z}(T)$ belongs to $\mathcal{C}^{\gam}(\red)$, for any 
$\gam\in\,]0,\alpha - \frac{d}{2}[$. 
\item 
If the hypotheses of Theorem \ref{t1} are satisfied with some $\alpha\in [0,\frac{d}{2}[$, then $E\big( \Vert
v_{G,Z}(T)\Vert^2_{L^p(\red)}\big) < \infty$,
for any $p\in\,]2,\frac{2d}{d-2\alpha}[$. 
\end{enumerate} 
\end{corol} 

\begin{obs} 
\label{o1} 
\begin{description} 
\item{\rm (a)} 
Let $G$ be as in Example \ref{e1}. Assume that condition (\ref{1.5}) holds for some $\alpha \in\,]\frac{d}{2}, k[$. 
Then the conclusion of part 1 of Corollary \ref{c1} holds. 
This applies for instance to the wave equation in dimension $d=1$. 

If condition (\ref{1.5}) holds for some $\alpha \in [0,\frac{d}{2}\wedge k[$, then the conclusion 
of part 2 of Corollary \ref{c1} holds. 
\item{\rm (b)} Let $G$ be as in part (a) and $\Gam (dx)$ be as in Example \ref{e2}. Suppose that 
$\frac{d}{2} < k-\frac{\beta}{2}$. Then a.s., $v_{G,Z}(T)$ belongs to $\mathcal{C}^{\gam}(\red)$, 
for any $\gam\in\,]0,k-\frac{\beta+d}{2}[$. 
\end{description} 
\end{obs} 

\section {\bf Path properties in time of the stochastic integral} 

We are now interested in the behaviour in $t$ of the sample paths of the process 
$u_{G,Z} = (u_{G,Z}(t),\ t\in[0,T])$, where 
\beqn%\label{rde16a}
   u_{G,Z}(t) = \int_0^t \int_{\red}G(t-s, \cdot-y) Z(s,y) M(ds,dy). 
\eeqn
We notice that in Theorem \ref{t1}, one can replace everywhere the finite time 
horizon $T$ by an arbitrary $t\in[0,T]$ and $G(s)$ by $G(t-s)$; 
therefore, under the assumptions of this theorem, the process $u_{G,Z}$ takes its values 
in the Hilbert space $H_2^{\alpha}(\red)$. 

Our aim is to prove H\"{o}lder continuity of the sample paths. 
We shall apply a version of Kolmogorov's continuity condition; hence we are led to estimate 
$L^p$-moments of stochastic integrals with values in a Hilbert space by means of an extension of
Burkholder's inequality.  We therefore need to 
identify the increasing process associated with the martingale $(v_{G,Z}(t),\ t\in[0,T])$ of Theorem \ref{t1}. We devote
the first part of this 
section to this problem; the second part deals with the study of H\"{o}lder continuity. 

\subsection{\bf The increasing process} 

Following \cite{metivier}, we term the {\it Meyer process} or {\it first increasing process} of the 
$H_2^{\alpha}(\red)$-valued martingale $(v_{G,Z}(t),\ t\in[0,T])$ the unique real-valued, continuous, 
increasing process, denoted by $(\langle v_{G,Z}\rangle_t,\ t\in[0,T])$, 
such that $\Vert v_{G,Z}(t)\Vert_{H_2^{\alpha}(\red)}^2 - \langle v_{G,Z}\rangle_t$ is a real-valued martingale. 

\begin{prop} 
\label{p2} 
Assume that the hypotheses of Theorem \ref{t1} are satisfied. Then for any $t\in [0,T]$, 
\begin{align} 
\langle v_{G,Z}\rangle_t &= \int_0^tds \int_{\red} d\xi\, |\tf Z(s)(\xi)|^2\nonumber \\ 
& \quad \times\int_{\red}\mu(d\eta) (1 + |\xi-\eta|^2)^{\alpha} \big|\tf G(s)(\xi - \eta)\big|^2. \label{2.1} 
\end{align} 
\end{prop} 

\noindent{\it Proof}. Assume first that $G$ and $Z$ satisfy 
the assumptions {\bf (G1')}, {\bf (G2)} and 
{\bf (G3)} with $T$ replaced by $t$, that is, $I_{G,Z}^{0,t} < \infty$, for any $t\in [0,T]$, where 
\beqn 
I_{G,Z}^{0,t} = \int_0^t ds \int_{\red} d\xi\, E(|\tf Z(s)(\xi)|^2) \int_{\red} \mu(d\eta)|\tf G(s) (\xi - \eta)|^2. 
\eeqn 
Suppose also that 
$G^{\alpha}(s) = (I - \Delta)^{\frac{\alpha}{2}}G(s)$ 
satisfies {\bf (G1')} and $I_{G,Z}^{\alpha,t}<\infty$, where 
\begin{align*} 
I_{G,Z}^{\alpha,t}& = \int_0^t ds \int_{\red} d\xi\, E(|\tf Z(s)(\xi)|^2)\\ 
&\quad \times \int_{\red} \mu(d\eta) (1 + |\xi - \eta|^2)^{\alpha}|\tf G(t-s) (\xi - \eta)|^2. 
\end{align*} 
Then, following Lemma 1 in \cite{dalangmueller}, 
for any $t\in[0,T]$ and $x\in \red$, the stochastic integral 
\beqn 
v_{(I-\Delta)^{\frac{\alpha}{2}}G,Z}(t,x) = 
\int_0^t\int_{\red} (I-\Delta)^{\frac{\alpha}{2}}G(s, x - y)(x) 
Z(s,y) M(ds,dy) 
\eeqn 
is well-defined as a Walsh stochastic integral. Its increasing process is given by 
\beqn 
\langle v_{(I-\Delta)^{\frac{\alpha}{2}}G,Z}\rangle _t = \int_0^tds \int_{\red} 
\mu(d\eta)\, | \tf \big( (I - \Delta)^{\frac{\alpha}{2}}G(s,x-\cdot) Z(s,\cdot)\big)(\eta)|^2, 
\eeqn 
(use Theorem 2.5 in \cite{walsh} and elementary properties of convolution and the Fourier transform). 

In particular, the process 
\beq 
\label{2.2} 
|v_{(I-\Delta)^{\frac{\alpha}{2}}G,Z}(t,x)|^2 - \int_0^tds \int_{\red} 
\mu(d\eta)\,\Big| \tf \big((I - \Delta)^{\frac{\alpha}{2}}G(s,x-\cdot) Z(s,\cdot)\big)(\eta)\Big|^2, 
\eeq 
$t\in[0,T]$, is a real-valued martingale. 

The properties of the Fourier transform yield 
\begin{align*} 
&\tf \big( (I - \Delta)^{\frac{\alpha}{2}}G(s,x-\cdot) Z(s,\cdot)\big)(\eta)\\ 
&\qquad = \Big ( \tf \big((I - \Delta)^{\frac{\alpha}{2}} G(s,x-\cdot)\big) * \tf Z(s,\cdot) \Big) (\eta)\\ 
&\qquad  = \tf \Big( (1 + |\xi'|^2)^{\frac{\alpha}{2}} \tf G(s,\cdot)(\eta-\xi ') \tf Z(s,\cdot)(\xi ')\Big) (x). 
\end{align*} 
Then, by Plancherel's theorem, 
\begin{align} 
& \int_{\red} dx \int_0^tds \int_{\red} 
\mu(d\eta)\, \Big| \tf \big((I - \Delta)^{\frac{\alpha}{2}}G(s,x-\cdot) Z(s,\cdot)\big)(\eta)\Big|^2\nonumber\\ 
&= \int_0^tds \int_{\red} d\xi\, |\tf Z(s)(\xi)|^2 
\int_{\red}\mu(d\eta) (1 + |\xi-\eta|^2)^{\alpha} \big|\tf G(s)(\xi - \eta)\big|^2.\label{2.3} 
\end{align} 

   Following Step 1 in the proof of Theorem \ref{t1}, 
\beqn 
\Vert v_{G,Z}(t)\Vert_{H_2^{\alpha}(\red)}^2 = \Vert v_{(I-\Delta)^{\frac{\alpha}{2}}G,Z}(t)\Vert^2_{L^2(\red)}, 
\eeqn 
for any $t\in[0,T]$. 

 Integrating over $x \in \red$ the expression in (\ref{2.2}) and using (\ref{2.3}), we find that the process 
\begin{align*} 
&\Vert v_{G,Z}(t)\Vert_{H_2^{\alpha}(\red)}^2 - 
\int_0^tds \int_{\red} d\xi\, |\tf Z(s)(\xi)|^2\nonumber \\ 
& \qquad \times\int_{\red}\mu(d\eta) (1 + |\xi-\eta|^2)^{\alpha} \big|\tf G(t-s)(\xi - \eta)\big|^2
\end{align*} 
is a real-valued martingale. This proves (\ref{2.1}) under the particular set of assumptions stated at the beginning of the
proof. 

   Assume next the setting of Step 2 in the proof of Theorem \ref{t1}. That is, $G(s)$ and $G^{\alpha}(s)$ 
satisfy {\bf(G1')} and condition (iii) holds. There exists a sequence of processes 
$(Z_n, n\ge 1)$ satisfying {\bf(G2)} and $I_{G,Z}^{\alpha,t} < \infty$, such that for any $t\in[0,T]$, 
\beq
\lim_{n\to\infty}E\big(\Vert v_{G,Z-Z_n}(t)\Vert^2_{H_2^{\alpha}(\red)}\big)
 = \lim_{n\to\infty}I_{G,Z_n-Z}^{\alpha,t} = 0.
 \label{2.4} 
\eeq
By the previous step, 
\begin{align*} 
M_t^n:&=\Vert v_{G,Z_n}(t)\Vert_{H_2^{\alpha}(\red)}^2 - 
\int_0^tds \int_{\red} d\xi\, |\tf Z_n(s)(\xi)|^2\nonumber \\ 
& \quad \times\int_{\red}\mu(d\eta) (1 + |\xi-\eta|^2)^{\alpha} \big|\tf G(t-s)(\xi - \eta)\big|^2, 
\end{align*} 
$t\in[0,T]$, is a real-valued martingale. Set 
\begin{align*} 
M_t& = \Vert v_{G,Z}(t)\Vert_{H_2^{\alpha}(\red)}^2 - 
\int_0^tds \int_{\red} d\xi\, |\tf Z(s)(\xi)|^2\nonumber \\ 
& \quad \times\int_{\red}\mu(d\eta) (1 + |\xi-\eta|^2)^{\alpha} \big|\tf G(t-s)(\xi - \eta)\big|^2. 
\end{align*} 

From (\ref{2.4}), it follows that 
$L^1(\Omega)$-$\lim_{n\to\infty} M_t^n = M_t$. This shows that 
$(M_t, t\in[0,T])$ is a martingale and proves (\ref{2.1}) in the setting of Step 2. 

   Finally, we consider the situation given by the hypotheses of the theorem.
 From Step 3 in the proof of Theorem \ref{t1}, it follows that for any $t\in[0,T]$, 
\beqn 
 \lim_{n\to\infty}E\big(\Vert v_{G-G_n,Z}(t)\Vert^2_{H_2^{\alpha}(\red)}\big) 
  = \lim_{n\to\infty}I_{G-G_n,Z}^{\alpha,t} = 0, 
\eeqn 
where $G_n(s) = G(s) * \psi_n$ and $(\psi_n, n\ge 1)$ is an approximation of the identity. 
The sequence $(G_n(s),\ n\ge 1)$ satisfies the conditions of the previous step. 
Therefore, we can conclude using a limiting procedure, in a manner analoguous to the previous step. This completes the
proof of the Proposition. 
\hfill \qed 

\begin{prop}
Assume the hypotheses of Theorem \ref{t1}. Fix $q \in [1,\infty[$. Then there is $C>0$ such that for all $t>0$,
\begin{align*}
   E\big(\Vert v_{G,Z}(t)\Vert_{H_2^{\alpha}(\red)}^{2q}\big)
   & \le C\, t^{q-1}
   \int_{0}^{t} ds\, E\left(\Vert Z(s)\Vert^{2q}_{L^2(\red)}\right)\\
  &\qquad\times \Big(\sup_{\xi\in\red} \int_{\red} \mu(d\eta)(1 + |\xi - \eta|^2)^{\alpha}|\tf G(s)(\xi - \eta)|^2\Big)^q.
\end{align*}
\label{rdp3}
\end{prop}

\noindent{\em Proof.} Using the Hilbert space version of Burkholder's inequality (\cite{metivier}, p.~212),
 Proposition \ref{p2}, H\"{o}lder's inequality and Plancherel's identity, we obtain 
\begin{align*} 
E\big(\Vert v_{G,Z}(t)\Vert_{H_2^{\alpha}(\red)}^{2q}\big) 
& \le C E\Big(\big(\int_{0}^{t} ds \int_{\red} d\xi\, |\tf Z(s)(\xi)|^2 
\int_{\red} \mu(d\eta)(1 + |\xi - \eta|^2)^{\alpha}\\ 
&\qquad\qquad\qquad\qquad\times |\tf G(s)(\xi - \eta)|^2\big)^q\Big)\\ 
&\le C t^{q-1} 
\int_{0}^{t} ds\, E\Big(\big(\int_{\red} d\xi\, |\tf Z(s)(\xi)|^2\\ 
&\qquad\qquad\times\int_{\red} \mu(d\eta)(1 + |\xi - \eta|^2)^{\alpha}|\tf G(s)(\xi - \eta)|^2\big)^q\Big)\\ 
&\le C t^{q-1} 
\int_{0}^{t} ds\, E\left(\Vert Z(s)\Vert^{2q}_{L^2(\red)}\right)\\
   &\qquad \times\Big(\sup_{\xi\in\red} \int_{\red} \mu(d\eta)(1 + |\xi - \eta|^2)^{\alpha}|\tf G(s)(\xi - \eta)|^2\Big)^q.
\end{align*} 
This proves the proposition.
\hfill \qed

\subsection{H\"{o}lder continuity in time} 

In this section, we consider the distribution-valued function $G(s)$ of Example 
\ref{e1}. Our goal is to give sufficient conditions ensuring a.s.-H\"{o}lder
 continuity of the sample paths of the process $(u_{G,Z}(t), t\in[0,T])$ defined in (\ref{0.3.1}).

We first study the case of a general covariance measure $\Gam$. In a second part, we consider the particular case
$\Gam(dx) = |x|^{-\beta}$, $\beta\in\,]0,d[$. The radial structure of 
this measure makes it possible to obtain a higher order of H\"{o}lder continuity. Indeed, we prove that the 
result obtained in this situation is optimal. 

\begin{teor} 
\label{t2.1} 
Let $\mathcal{L} = \partial_{tt}^2 + (-\Delta)^k$, $k\in\,]0,\infty[$, and let $G$ be the fundamental solution of
$\mathcal{L}u = 0$. Fix $\alpha \in [0, k[$ and assume that there exists $\eta\in\,]\frac{\alpha}{k},1[$ such that 
\beq 
\label{2.5.0} 
\int_{\red}\frac{\mu (d\xi)}{(1 + |\xi|^2)^{k\eta - \alpha}} < \infty. 
\eeq 
Fix $q\in [2, \infty[$ and assume that $\sup_{0\le s\le T}E(\Vert Z(s)\Vert^q_{L^2(\red)}) < \infty$. 
Then the $H_2^{\alpha}(\red)$-valued stochastic integral process $(u_{G,Z}(t),\  
0\le t\le T)$ satisfies 
\beq 
\label{2.5} 
E\big(\Vert u_{G,Z}(t_2) - u_{G,Z}(t_1)\Vert_{H_2^{\alpha}(\red)}^{q}\big) \le C (t_2 - 
t_1)^{q\big(\frac{1}{2}\wedge(1-\eta)\big)}, 
\eeq 
for any $0\le t_1\le t_2\le T$. Consequently, $(u_{G,Z}(t), 0\le t\le T)$ 
is $\gam$-H\"{o}lder continuous, for each $\gam\in\,\big]0, (\frac{1}{2}\wedge(1-\eta))-\frac{1}{q}\big[$. 
\end{teor} 

\noindent{\it Proof}. 
Fix $0\le t_1\le t_2\le T$ and set $q=2p$. Then 
\beqn 
E\big(\Vert u_{G,Z}(t_2) - u_{G,Z}(t_1)\Vert_{H_2^{\alpha}(\red)}^{2p}\big) \le C\big(T_1(t_1,t_2) + T_2(t_1,t_2)\big),
\eeqn 
where 
\begin{align*} 
T_1(t_1,t_2) & = E \big(\Vert \int_{t_1}^{t_2} G(t_2 - s, \cdot - y) Z(s,y) M(ds,dy)
\Vert^{2p}_{H_2^{\alpha}(\red)}\big),\\ 
T_2(t_1,t_2) & = E \big(\Vert \int_0^{t_1} \big(G(t_2-s, \cdot - y) - G(t_1-s, \cdot - y)\big)\\ 
& \qquad\qquad \times Z(s,y) M(ds,dy)\Vert^{2p}_{H_2^{\alpha}(\red)}\big). 
\end{align*} 
Arguing as in Proposition \ref{rdp3} and using (\ref{1.6.1}), we obtain 
\begin{align} 
T_1(t_1,t_2) 
%& \le C E\Big(\Big(\int_{t_1}^{t_2} ds \int_{\red} d\xi\, |\tf Z(s)(\xi)|^2 
%\int_{\red} \mu(d\eta)(1 + |\xi - \eta|^2)^{\alpha}\nonumber\\ 
%&\quad\times |\tf G(t_2-s)(\xi - \eta)|^2\Big)^p\Big)\nonumber\\ 
&\le C (t_2 - t_1)^{p-1} 
\int_{t_1}^{t_2} ds\, E\Big(\Vert Z(s)\Vert^{2p}_{L^2(\red)}\Big) \nonumber\\ 
&\qquad\times\left(\sup_{\xi\in\red}\int_{\red} \mu(d\eta)(1 + |\xi - \eta|^2)^{\alpha}|\tf G(t_2-s)(\xi -
\eta)|^2\right)^p \nonumber\\ 
&\le C (t_2 - t_1)^p 
\sup_{0\le s\le T}E\big(\Vert Z(s)\Vert^{2p}_{L^2(\red)}\big)\nonumber\\ 
&\qquad\times \sup_{0\le s\le T} \sup_{\xi\in\red}\left( \int_{\red} \mu(d\eta) (1 + |\xi - \eta|^2)^{\alpha}|\tf G(s)(\xi
- 
\eta)|^2 \right)^p\nonumber\\ 
&\le C (t_2 - t_1)^p. \label{2.6} 
\end{align} 

   We now study the contribution of $T_2(t_1,t_2)$. Clearly, 
\beqn 
\int_{\red} \mu(d\eta) (1 + |\xi - \eta|^2)^{\alpha}|\tf (G(t_2 -s) - G(t_1 - s)) (\xi - \eta)|^2 \le I_1(t_1,t_2) + 
I_2(t_1,t_2), 
\eeqn 
where 
\begin{align*} 
I_1(t_1,t_2)& = \int_{|\xi - \eta |\le 1}\mu(d\eta) (1 + |\xi - \eta|^2)^{\alpha} 
|\tf (G(t_2 -s) - G(t_1 - s))(\xi - \eta)|^2 ,\\ 
I_2(t_1,t_2)& = \int_{|\xi - \eta | > 1}\mu(d\eta) (1 + |\xi - \eta|^2)^{\alpha}|\tf (G(t_2 -s) - G(t_1 - s))(\xi - \eta)|^2. 
\end{align*} 
By (\ref{rde11a}), the mean-value theorem, the bound (\ref{1.6}) and assumption (\ref{2.5.0}), 
\begin{align*} 
I_1(t_1,t_2)&\le (t_2 - t_1)^2 \int_{|\xi - \eta |\le 1} \mu(d\eta) (1 + |\xi - \eta|^2)^{\alpha}\\ 
%& \le 2^k (t_2 - t_1)^2 \int_{\red}\frac{\mu(d\eta)}{(1 + |\xi - \eta|^2)^k}\\ 
& \le 2^k (t_2 - t_1)^2 \sup_{\xi\in\red} \int_{\red}\frac{\mu(d\eta)}{(1 + |\xi - \eta|^2)^{k-\alpha}}\\ 
& \le C (t_2 - t_1)^2 \int_{\red}\frac{\mu(d\eta)}{(1 + |\eta|^2)^{k-\alpha}}\\ 
& \le C (t_2 - t_1)^2. 
\end{align*} 
By the formula $\sin x - \sin y = 2 \cos \frac{x+y}{2} \sin\frac{x-y}{2}$, 
\begin{align*} 
I_2(t_1,t_2)& \le \int_{|\xi - \eta | > 1}\mu(d\eta) (1 + |\xi - \eta|^2)^{\alpha} 
\frac{\Big (\sin\big(\frac{1}{2} (t_2 - t_1) |\xi - \eta |^k\big)\Big)^{2(1-\eta)}}{|\xi - \eta|^{2k}}\\ 
&\le C (t_2 - t_1)^{2(1-\eta)} \int_{\red} \frac{\mu(d\eta)}{(1 + |\xi - \eta|^2)^{k\eta-\alpha}}\\ 
&\le C (t_2 - t_1)^{2(1-\eta)}. 
\end{align*} 
Consequently, 
\begin{align} 
&\sup_{0\le s\le T} \sup_{\xi\in\red}\int_{\red} \mu(d\eta) 
(1 + |\xi - \eta|^2)^{\alpha}|\tf (G(t_2 -s) - G(t_1 - s)) (\xi - 
\eta)|^2\nonumber\\ 
&\le C (t_2 - t_1)^{2(1-\eta)}\label{2.7}. 
\end{align} 

Using (\ref{2.7}) and arguing as in the lines that led to (\ref{2.6}), we see that
\begin{align} 
T_2(t_1,t_2)
%& \le C \int_0^{t_1} ds\, E\Big( \int_{\red} d \xi\, |\tf Z(s)(\xi)|^2 
%\int_{\red} \mu(d\eta) (1 + |\xi - \eta|^2)^{\alpha}\nonumber\\ 
%&\quad \times|\tf (G(t_2 -s) - G(t_1 - s))(\xi - \eta)|^2\Big)^p\nonumber\\ 
&\le C \sup_{0\le s\le T} E(\Vert Z(s)\Vert^{2p}_{L^2(\red)}) 
\sup_{0\le s\le T}\sup_{\xi\in\red} 
\Big(\int_{\red} \mu(d\eta) (1 + |\xi - \eta|^2)^{\alpha}\nonumber\\ 
& \qquad\times|\tf (G(t_2 -s) - G(t_1 - s))(\xi - \eta)|^2\Big)^p\nonumber\\ 
&\le C (t_2 - t_1)^{2p(1-\eta)}\label{2.8}. 
\end{align} 
Finally, (\ref{2.5}) is a consequence of (\ref{2.6}) and (\ref{2.8}). 

The statement on H\"{o}lder continuity follows from 
Kolmogorov's continuity condition \cite[Chap.I, \S 2]{RY}. \hfill \qed 
\medskip 

The previous theorem, together with part 1 of Corollary \ref{c1}, yields the following. 

\begin{corol} 
\label{c2} 
Suppose that the hypotheses of Theorem \ref{t2.1} are satisfied with some 
$\alpha \in \,]\frac{d}{2}, \infty[$ and $\eta\in\,]\frac{\alpha}{k},1[$ Then there is a version of the process 
\beqn 
\big(u_{G,Z}(t,x),\ (t,x)\in[0,T]\times \red\big) 
\eeqn 
that belongs to 
$\mathcal{C}^{\gam_1,\gam_2}([0,T]\times\red)$, with $\gam_1 \in \,]0, (\frac{1}{2}\wedge(1-\eta))-\frac{1}{q}[$ and
$\gam_2\in \,]0,\alpha-\frac{d}{2}[$. 
\end{corol} 
\medskip 

Consider now the particular case $\Gam(dx) = |x|^{-\beta}$, $\beta\in\,]0,d[$. The results obtained in Theorem \ref{t2.1}
can be improved as follows. 

\begin{teor} 
\label{t2.2} 
Let $\mathcal{L}$ and $G$ be as in Theorem \ref{t2.1}. Fix $\alpha \in [0,\infty[$, $k>\alpha$ and assume that 
$\Gam(dx) = |x|^{-\beta}$ with $\beta\in \,]0, 2(k - \alpha)[$. Fix $q\in[2,\infty[$ and suppose that $\sup_{0\le s\le
T}E(\Vert Z(s)\Vert^q_{L^2(\red)}) < \infty$. 
Then for any $0\le t_1 \le t_2 \le 1$, 
\beq 
\label{2.9} 
E\big(\Vert u_{G,Z}(t_2) - u_{G,Z}(t_1)\Vert^q_{H_2^{\alpha}}\big) \le C (t_2-t_1)^{q(1-\frac{\beta+2\alpha}{2k})}. 
\eeq 
Consequently, $(u_{G,Z}(t),\ 0\le t\le T)$ is $\gam$-H\"{o}lder continuous for any 
$\gam \in \,\big]0, (1-\frac{\beta+2\alpha}{2k})-\frac{1}{q}\big[$. 
\end{teor} 

\begin{remark} If $\beta+2\alpha < k$, then we obtain a stronger conclusion than in Theorem \ref{t2.1}. 
\end{remark}

\bigskip 

\noindent{\it Proof}. As in the proof of Theorem \ref{t2.1}, set $q=2p$ and 
\begin{align*} 
T_1(t_1, t_2) & = E \big(\Vert \int_{t_1}^{t_2} G(t_2 - s, \cdot - y) Z(s,y) M(ds,dy)
\Vert^{2p}_{H_2^{\alpha}(\red)}\big),\\ 
T_2(t_1, t_2) & = E \big(\Vert \int_0^{t_1} \big(G(t_2-s, \cdot - y) - G(t_1-s, \cdot - y)\big) Z(s,y) 
M(ds,dy)\Vert^{2p}_{H_2^{\alpha}(\red)}\big). 
\end{align*} 
Let 
\beqn 
T_{11}(t_1, t_2) = \int_0^{t_2-t_1} ds \sup_{\xi\in\red} \big( \int_{\red} \mu(d\eta) 
(1 + |\xi - \eta|^2)^{\alpha}|\tf G(s)(\xi - \eta)|^2\big) ^p. 
\eeqn 
Then, proceeding as in the steps that led to (\ref{2.6}), we find that
\beqn 
T_1(t_1, t_2) \le C(t_2 - t_1)^{p-1}\sup_{0\le s \le T} E\big(\Vert Z(s)\Vert^{2p}_{L^2(\red)}\big) 
T_{11}(t_1,t_2). 
\eeqn 
Introducing the new variables $(\tilde \xi, \tilde \eta) = s^{\frac{1}{k}}(\xi, \eta)$ and 
substituting $|\eta|^{-d+\beta}$ for $\mu (d\eta)$ and formula (\ref{rde11a}) for $\tf G(s)$  yields 
\begin{align*} 
T_{11}(t_1, t_2) & = \int_0^{t_2-t_1} ds\, s^{p(2-\frac{\beta + 2\alpha}{k})}\\ 
&\quad\times \sup_{\tilde \xi \in \red}\Big( \int_{\red} \frac{d \tilde\eta}{|\tilde\eta|^{d-\beta}} \,
(s^{\frac{2}{k}} + |\tilde \xi - \tilde \eta|^2)^{\alpha}\frac{\sin^2( |\tilde \xi - \tilde \eta|^k)} 
{ |\tilde \xi - \tilde \eta|^{2k}}\Big)^p. 
\end{align*} 
Taking into account (\ref{1.5.1}), (\ref{1.6}) for $\mu(d\eta) = |\eta|^{-d+\beta} d\eta$ and the 
remark made in Example \ref{e2}, we obtain 
\beqn 
T_{11}(t_1, t_2) \le C (t_2-t_1)^{p(2-\frac{\beta+2\alpha}{k})+1}. 
\eeqn 
Consequently, 
\beq 
\label{2.10} 
T_1(t_1, t_2) \le C(t_2 - t_1)^{p(3-\frac{\beta+2\alpha}{k})}. 
\eeq 
For the analysis of the term $T_2(t_1, t_2)$, we also follow the same scheme as in the 
proof of Theorem \ref{t2.1} but we improve the upper bound on $I_2(t_1, t_2)$, as follows. 
Set $h=t_2-t_1$ and consider the change of variables 
$(\tilde \xi, \tilde \eta) = (\frac{h}{2})^{\frac{1}{k}}(\xi, \eta)$. 
Then, 
\begin{align*} 
I_2(t_1, t_2)& \le \int_{|\xi -\eta|>1}\, \frac{d\eta}{|\eta|^{d-\beta}}\big(1+|\xi-\eta|^2\big)^{\alpha} 
\frac{\sin^2 \big(\frac{1}{2}(t_2-t_1)|\xi-\eta|^k\big)}{|\xi-\eta|^{2k}}\\ 
& \le C h^{2-\frac{\beta+2\alpha}{k}}\int_{|\tilde\xi-\tilde\eta|>(\frac{h}{2})^{\frac{1}{k}}} 
\frac{d\tilde\eta}{|\tilde \eta|^{d-\beta}}(1+|\tilde \xi - \tilde \eta|^2)^{\alpha} 
\frac{\sin^2(|\tilde \xi-\tilde\eta|^k)}{|\tilde\xi - \tilde \eta|^{2k}}, 
\end{align*} 
where $C$ is a constant depending on $k$. Therefore, 
\beq 
\label{2.11} 
T_2(t_1,t_2) \le C (t_2 - t_1)^{p(2-\frac{\beta+2\alpha}{k})}. 
\eeq 
The estimates (\ref{2.10}) and (\ref{2.11}) imply (\ref{2.9}). The proof of the theorem is complete. 
\hfill \qed 
\bigskip 

We finish this section by showing that Theorem \ref{t2.2} provides an optimal result. We do this by studying the case
where $Z$ is the smooth deterministic function $Z(s,x) = e^{-|x|^2/2}$, with no dependence on $s$. For this $Z$, we shall
write $u_G(t)$ instead of $u_{G,Z}(t)$. 

\begin{teor} 
\label{t2.3} 
Let $\mathcal{L}$, $G$ and $\Gam$ be as in Theorem \ref{t2.2}. Fix $t_0\in (0,1]$ and assume $\beta
\in\,]0,2(k-\alpha)[$. Then there exists a constant $C>0$ such that for any 
$t_1, t_2$ satisfying $t_0 \le t_1 \le t_2\le 1$,  
\beq 
\label{2.12} 
E\big(\Vert u_G(t_2) - u_G(t_1)\Vert^2_{H_2^{\alpha}(\red)}\big) \ge C |t_2 - t_1|^{2-\frac{\beta+2\alpha}{k}}. 
\eeq 
Consequently, a.s. the mapping $t\mapsto u_G(t)$ is {\em not} $\gam$-H\"{o}lder continuous for 
$\gam > 1-(\beta+2\alpha)/(2k)$, while it is $\gam$-H\"{o}lder continuous for 
$\gam <1-(\beta+2\alpha)/(2k)$. 
\end{teor} 

\noindent{\it Proof}. Let $p(\xi)$ denote the standard Gaussian density function. Using the isometry property (\ref{1.1}),
we obtain 
\beqn 
E\big(\Vert u_G(t_2) - u_G(t_1)\Vert^2_{H_2^{\alpha}(\red)}\big) \ge S(t_1, t_2), 
\eeqn 
where 
\begin{align} 
\nonumber
S(t_1, t_2) &= \int_0^{t_1} ds\int_{\red} d\xi\, p(\xi)^2 \int_{\red} \mu(d\eta) \big(1+|\xi - \eta|^2\big)^{\alpha} \\
   &\qquad\qquad\qquad \times\big|\tf G(t_2-s)(\xi - \eta) - \tf G(t_1-s)(\xi - \eta)\big|^2. 
\label{2.13} 
\end{align} 
Set $h=(t_2-t_1)/2$. 
By the formula $\sin x - \sin y = 2 \cos \frac{x+y}{2} \sin \frac{x-y}{2}$, Fubini's 
theorem and integrating with respect to the time variable $s$, we obtain 
\begin{align*} 
S(t_1, t_2) &= 4 \int_{\red} d\xi\, p(\xi)^2 \int_{\red} \frac{d\eta}{|\eta|^{d-\beta}} \big(1+|\xi - \eta|^2\big)^{\alpha} 
\frac{\sin^2(h|\xi - \eta|^k)}{|\xi - \eta|^{2k}}\\ 
&\qquad\qquad\times \left( \frac{t_1}{2} - \frac{\sin ((t_2-t_1)|\xi - \eta|^k)}{4|\xi - \eta|^k} 
+ \frac{\sin((t_1+t_2)|\xi - \eta|^k)}{4|\xi - \eta|^k}\right). 
\end{align*} 
Notice that for $|\xi - \eta|^k > 2/t_1$ and, in particular, for $|\eta| > 2 (2/t_1)^{1/k}$ and $|\xi| < (2/t_1)^{1/k}$, the
factor in parentheses is bounded below by 
$t_1/4$. Therefore,
\begin{align*} 
   S(t_1, t_2)&\geq t_0 \int_{|\xi| < (2/t_1)^{1/k}} d\xi\, p(\xi)^2 \int_{|\eta| > 2 (2/t_1)^{1/k}}
\frac{d\eta}{|\eta|^{d-\beta}}\\
   &\qquad\qquad\qquad\qquad\times \big(1+|\xi - \eta|^2\big)^{\alpha}\, \frac{\sin^2(h|\xi - \eta|^k)}{|\xi - \eta|^{2k}}.
\end{align*} 
Let $a = (2/t_1)^{1/k}$. Note that 
\beqn
   \{(\xi,\eta): |\xi| < a,\ |\eta| > 2 a\} \supset
   \{(\xi,\eta): |\xi| < a,\ |\xi -\eta| > 3 a\}
\eeqn
and that $1/|\eta| > 1/(2|\xi -\eta|)$ for $(\xi,\eta)$ in these sets. With this inequality and this smaller domain of
integration, we use the change of variables  $\tilde \eta = h^{\frac{1}{k}} (\xi - \eta)$ ($\xi$ fixed) to see that 
\begin{align*} 
S(t_1, t_2) &\ge \frac{t_0}{2} h^{2-\frac{\beta+2\alpha}{k}} \int_{|\xi| < a} d\xi\, p(\xi)^2  \int_{|\tilde \eta|^k >  3^k
a^k} 
\frac{d\tilde \eta}{|\tilde\eta|^{d-\beta+2k}}(h^{\frac{2}{k}}+|\tilde\eta|^2)^{\alpha} 
\sin^2(|\tilde\eta|)\\ 
&\ge \frac{t_0}{2} h^{2-\frac{\beta+2\alpha}{k}} \int_{|\xi| < 2^{1/k}} d\xi\, p(\xi)^2\int_{|\tilde\eta|^k >3^k 2/t_0} 
\frac{d\tilde\eta}{|\tilde\eta|^{d-\beta +2(k-\alpha)}} \sin^2(|\tilde\eta|). 
\end{align*} 
Notice that the last double integral is a positive finite constant. Hence, the inequality (\ref{2.12}) is proved. 

   We now use the fact that $u_G$ is a Gaussian stationary process together with classical results
 on Gaussian processes to translate the lower bound (\ref{2.12}) into a statement concerning absence of
H\"older continuity of the sample paths of $t \mapsto u_G(t)$. Fix $\gamma \in \,]1-(\beta+2\alpha)/(2k), 1]$ and
assume by
contradiction that for almost all $\omega$, there is $C(\omega) < \infty$ such that for all $t_0 \leq t_1
< t_2 \leq 1$, 
\begin{align*} 
  & \sup_{t_1 < t_2}\ \sup_{\varphi \in H_2^{-\alpha}(\red),\, \varphi \not\equiv 0} \frac{\langle u_G(t_2) - u_G(t_1),\
\varphi \rangle}{(t_2 - t_1)^\gamma\, \Vert \varphi \Vert_{H_2^{-\alpha}(\red)} } \\
  &\qquad\qquad = \sup_{t_1 < t_2}\ \frac{\Vert u_G(t_2) - u_G(t_1)\Vert_{H_2^{\alpha}(\red)}  }{(t_2
- t_1)^\gamma }\\
  &\qquad\qquad < C(\omega).
\end{align*} 
Then the real-valued Gaussian stochastic process
\beqn
   \left(\frac{\langle u_G(t_2) - u_G(t_1),\ \varphi \rangle}{(t_2 - t_1)^\gamma\, \Vert \varphi
\Vert_{H_2^{-\alpha}(\red)} },\quad t_1 < t_2,\ \varphi \in H_2^{-\alpha}(\red),\, \varphi \not\equiv 0\right) 
\eeqn
is finite a.s. By Theorem 3.2 of \cite{adler}, it follows that
\beqn
   E\left(\sup_{t_1 < t_2}\ \sup_{\varphi \in H_2^{-\alpha}(\red),\,
 \varphi \not\equiv 0} \left(\frac{\langle u_G(t_2) - u_G(t_1),\ \varphi \rangle}{(t_2 - t_1)^\gamma\,
\Vert \varphi \Vert_{H_2^{-\alpha}(\red)} } \right)^2 \right) < \infty.
\eeqn
Thus,
\beqn
   E\left(\sup_{t_1 < t_2}\ \frac{\Vert u_G(t_2) - u_G(t_1)\Vert^2_{H_2^{\alpha}(\red)}}{(t_2 -
t_1)^{2\gamma} }\right) < \infty.
\eeqn
In particular, there would exist $K< \infty$ such that
\beqn
  E\left( \Vert u_G(t_2) - u_G(t_1)\Vert^2_{H_2^{\alpha}(\red)} \right) \leq K\, |t_2 - t_1|^{2\gamma}.
\eeqn
However, this would contradict (\ref{2.12}) since $2\gamma > 2 - (\beta + 2 \alpha) / k$. We conclude that $t \mapsto
u_G(t)$ is {\em not} $\gamma$-H\"older continuous for $\gamma > 1-(\beta+2\alpha)/(2k)$.

   On the other hand, for $\gam < 1-(\beta+2\alpha)/(2k)$, the map $t \mapsto u_G(t)$ is $\gam$-H\"older continuous by
Theorem \ref{t2.2}, since in this theorem, $q$ can be taken arbitrarily large.
\hfill \qed

%%%%%%%%%%%%%%%%%%% 
%%%%%%%%%%%%%%%%%%%% 
\section {\bf Application to stochastic partial differential equations} 

This section is devoted to studying the properties of the sample paths of the solution of 
the spde 
\begin{align} 
&\left(\frac{\partial^2}{\partial t^2} + (- \Delta)^{(k)}\right) u(t,x) = \sigma\big(u(t,x)) \dot F(t,x) 
+ b \big(u(t,x)\big),\nonumber\\ 
&u(0,x) = v_0(x),\qquad 
\frac{\partial}{\partial t}u(0,x) = \tilde v_0(x). \label{3.1} 
\end{align} 
In this equation, $t \in [0,T]$ for some fixed $T >0$, and $x \in \red$. We assume that $k\in\,]0,\infty[$ ($k$ is not
necessarily an integer), $\sigma$ and $b$ are Lipschitz continuous functions and moreover, that 
\beq 
\label{3.2} 
|\sigma(z)| + |b(z)| \le C|z|, 
\eeq 
for some positive constant $C>0$. Notice that the assumption (\ref{3.2}) was also made in \cite{dalangmueller} and it is
also standard in the study of the deterministic wave equation (see for instance \cite[Chapter 6]{hormander} or
\cite{sogge}). 

Concerning the initial conditions, we assume for the moment that $v_0\in L^2(\red)$, $\tilde v_0 \in H_2^{-k}(\red)$.
Regarding the noise $\dot F$, we assume that its spectral measure satisfies (\ref{1.5}).

   By a solution of (\ref{3.1}), we mean an $L^2(\red)$-valued stochastic process $\big(u(t),\ 0\le t \le T \big)$ 
satisfying $\sup_{0\le t\le T} E\big(\Vert u(t)\Vert^2_{L^2(\red)}\big) < \infty$ and 
\begin{align} 
u(t,\cdot) &= \frac{d}{dt} G(t)*v_0 + G(t)*\tilde v_0\nonumber\\ 
&\qquad + \int_0^t \int_{\red} G(t-s,\cdot - y) \sigma (u(s,y)) M(ds,dy)\nonumber\\ 
&\qquad + \int_0^t ds \int_{\red} dy\, G(t-s,\cdot - y) b(u(s,y)). \label{3.3} 
\end{align} 
Here, $G$ is the fundamental solution of $\mathcal{L} f = 0$, where $\mathcal{L} = \big(\partial_{tt}^2 + (-
\Delta)^{(k)}\big)$, and the stochastic integral is of the type considered in the preceding sections (see also 
Section 2 in \cite{dalangmueller}). 

The path integral is also well-defined. Indeed, let $Z= (Z(s),\ s\in[0,T])$ be a stochastic processes satisfying the conditions
stated at the beginning of Section \ref{sec2}, and let $G: [0,T] \rightarrow \mathcal{S}'(\red)$ be such that for any
$s\in[0,T]$, 
$\tf G(s)$ is a function and 
\beq 
\label{3.4} 
\int_0^T ds \sup_{\xi\in\red} |\tf G(s)(\xi)|^2 < \infty. 
\eeq 
Then for any $t\in[0,T]$ a.s., 
\beqn 
   x \mapsto J_{G,Z}(t,x): = \int_0^t ds\, \big(G(s) * Z(s)\big) (x) 
\eeqn 
defines an $L^2(\red)$-valued function. Moreover, 
\beq 
\label{3.5} 
\Vert J_{G,Z}(t)\Vert^2_{L^2(\red)} \le C \int_0^t ds\, \Vert Z(s)\Vert^2_{L^2(\red)}\, \sup_{\xi\in\red}|\tf
G(s)(\xi)|^2 . 
\eeq 

Assume the following condition, which is stronger than (\ref{3.4}): 
\beq 
\label{3.6}\int_0^T ds \sup_{\xi\in\red} (1 + |\xi|^2)^{\alpha}|\tf G(s)(\xi)|^2 < \infty, 
\eeq 
for some $\alpha\in[0,\infty[$. 
Easy computations based on Fubini's theorem yield 
\beqn 
\Vert J_{G,Z}(t)\Vert^2_{H_2^{\alpha}(\red)} = \Vert J_{(I -
\Delta)^{\frac{\alpha}{2}}G,Z}(t)\Vert^2_{L^2(\red)}\qquad\mbox{a.s.} 
\eeqn 

Fix $q\in[2,\infty[$. Schwarz's inequality and Fubini's theorem yield a.s.: 
\begin{align*} 
\Vert J_{(I-\Delta)^{\frac{\alpha}{2}}G,Z}(t)\Vert^q_{L^2(\red)} &\le C \left( \int_0^t ds\,  
\Vert (I-\Delta)^{\frac{\alpha}{2}}G(s) 
* Z(s)\Vert^2_ {L^2(\red)}\right)^{\frac{q}{2}}\nonumber\\ 
& =C \left( \int_0^t ds \int_{\red}d\xi\, (1+|\xi|^2)^\alpha |\tf G(s)(\xi)|^2 |\tf Z(s)(\xi)|^2
\right)^{\frac{q}{2}}\nonumber\\ 
&\le C \left( \int_0^t ds\, \Vert Z(s)\Vert^2_ {L^2(\red)}\, \sup_{\xi\in\red}\big((1+|\xi|^2)^{\alpha}|\tf
G(s)(\xi)|^2\big) \right)^{\frac{q}{2}}\nonumber\\ 
&\le C \int_0^t ds\, \Vert Z(s)\Vert^q_{L^2(\red)}\, \sup_{\xi\in\red}\big((1+|\xi|^2)^{\frac{\alpha q}{2}}|\tf
G(s)(\xi)|^q\big) 
. 
\end{align*} 

   Now let $G$ be the fundamental solution of $\mathcal{L} f= 0$. 
Assume that $\alpha\in [0,k[$ and $\sup_{t\in[0,T]} E(\Vert Z(t)\Vert^q_{L^2(\red)}) < \infty$. 
Then
\beqn
   \sup_{\xi\in\red}(1+|\xi|^2)^{\frac{\alpha q}{2}}|\tf G(s)(\xi)|^q 
      \leq C \sup_{\xi\in\red}(1+|\xi|^2)^{q(\alpha - k)/2} < \infty,
\eeqn
and therefore the above inequalities yield 
\begin{align} 
E\big(\Vert J_{(I-\Delta)^{\frac{\alpha}{2}}G,Z}(t)\Vert^q_{L^2(\red)}\big) &\le C 
\int_0^t ds \sup_{\xi\in\red}(1+|\xi|^2)^{\frac{\alpha q}{2}}|\tf G(s)(\xi)|^q\nonumber\\ 
&\qquad\qquad\times E(\Vert Z(s)\Vert^q_{L^2(\red)})\nonumber\\ 
& < \infty. \label{3.8.1} 
\end{align} 

    Set 
\beqn 
J_{b}(t) = \int_0^t ds \int_{\red} dy\, G(t-s,\cdot - y) b(u(s,y)) 
= \int_0^t ds\, \big(G(s) * b(u(s)\big). 
\eeqn 
Particularizing (\ref{3.8.1}) to $\alpha=0$, $q=2$ and $Z(s,x) = b(u(t-s,x))$ yields 
\beqn 
%\label{3.8.2} 
E\big(\Vert J_{b}(t) \Vert^2_ {L^2(\red)}\big) \le 
C \sup_{0\le s\le T}E\big( \Vert u(s)\Vert^2_{L^2(\red)}\big) \int_0^T ds \sup_{\xi\in\red}|\tf 
G(s)(\xi)|^2 < \infty. 
\eeqn 

\medskip

   A slight extension of Theorem 9 in \cite{dalangmueller} provides the existence of a unique 
solution of equation (\ref{3.3}), in the sense given above. 
We observe that in \cite{dalangmueller}, $k\in{\mathbb{N}}$ and $b=0$. 

   By means of Burkholder's and H\"{o}lder's inequalities (as in the calculation that led to (\ref{2.6})), the inequality
(\ref{3.8.1}) with $\alpha = 0$ and a version of Gronwall's lemma (see \cite[Lemma 15]{dalang}), one can easily show that
for any $q\in[2,\infty[$, 
\beq 
\label{3.9} 
\sup_{0\le s\le T} E\big( \Vert u(s)\Vert^q_{L^2(\red)}\big) < \infty. 
\eeq 

   In the next theorem, we analyze the existence of $H^{\alpha}_2(\red)$-valued solutions to (\ref{3.3}). 

\begin{teor} 
\label{t3} 
Let $\sig$, $b$ be real-valued Lipschitz continuous functions satisfying (\ref{3.2}). 
Fix $\alpha\in [0,k[$ and assume that 
$v_0\in H^{\alpha}_2(\red)$, $\tilde v_0\in H_2^{\alpha - k}(\red)$ and 
\beqn 
\int_{\red} \frac{\mu (d\xi)}{(1 + |\xi|^2)^{k-\alpha}} < \infty. 
\eeqn 
Then for any $q\in [2,\infty[$, the solution of (\ref{3.3}) satisfies 
\beq 
\label{3.10} 
\sup_{0\le t\le T} E \big(\Vert u(t)\Vert^q_{H^{\alpha}_2(\red)}\big) < \infty. 
\eeq  
\end{teor} 

\noindent{\it Proof}. Fix $q\in [2,\infty[$. We shall check that each term on the right hand side of (\ref{3.3}) belongs to 
$L^q(\Omega; H^{\alpha}_2(\red))$, with norm uniformly bounded over $t\in[0,T]$. 

   Set $U_1(t) = \frac{d}{dt} G(t) * v_0$. Then 
\begin{align*} 
\Vert U_1(t)\Vert_{H^{\alpha}_2(\red)} 
&= \big \Vert (1+|\cdot|^2)^{\frac{\alpha}{2}}\tf \big( \frac{d}{dt} G(t)(\cdot)\big) \tf v_0(\cdot) \big \Vert_
{L^2(\red)}\\ 
& =\big \Vert (1+|\cdot|^2)^{\frac{\alpha}{2}} \cos(t|\cdot|^k)\tf v_0(\cdot) \big \Vert_ {L^2(\red)}\\ 
&\le \Vert v_0\Vert_{H^{\alpha}_2(\red)}. 
\end{align*} 
Similarly, define $U_2(t) = G(t) * \tilde v_0$. Then, by (\ref{rde12a}), 
\begin{align*} 
\Vert U_2(t)\Vert^2_{H^{\alpha}_2(\red)} 
& =\big \Vert (1+|\cdot|^2)^{\frac{\alpha}{2}} \tf G(t)(\cdot) \tf \tilde v_0(\cdot)\big \Vert^2_ {L^2(\red)}\\ 
& \le 2^k(1+T^2) \Vert \tilde v_0\Vert^2_{H^{\alpha-k}_2(\red)}. 
\end{align*} 
Hence, 
\beq 
\label{3.11} 
\sup_{0\le t \le T}  (\Vert U_1(t)\Vert^2_{H^{\alpha}_2(\red)} + \Vert U_2(t)\Vert^2_{H^{\alpha}_2(\red)}) < \infty. 
\eeq 

   Let 
\beqn
U_3(t) = \int_0^t \int_{\red} G(t-s,\cdot - y) \sigma(u(s,y)) M(ds,dy).
\eeqn
Using (\ref{3.2}), we see as in (\ref{2.6}) that 
\begin{align*} 
E\big(\Vert U_3(t)\Vert^q_{H^{\alpha}_2(\red)}\big)& \le C \sup_{0\le s\le T} E(\Vert u(s)\Vert^q_{L^2(\red)})\\ 
&\quad \times \sup_{0\le s\le T} \sup_{\xi\in\red}\Big(\int_{\red} \mu(d\eta) (1+|\xi - \eta|^2)^{\alpha} 
|\tf G(s)(\xi - \eta)|^2\Big)^{\frac{q}{2}}. 
\end{align*} 
By (\ref{1.6.1}) and (\ref{3.9}),  
\beq 
\label{3.12} 
\sup_{0\le t \le T} E\big(\Vert U_3(t)\Vert^q_{H^{\alpha}_2(\red)}\big) < \infty. 
\eeq 

   Finally, set 
\beqn U_4(t) = \int_0^t ds \int_{\red} dy\, G(t-s,\cdot - y) b(u(s,y)).
\eeqn 
The estimate (\ref{3.8.1}), (\ref{3.2}) and (\ref{3.9}) imply 
\beq 
\label{3.13} 
\sup_{0\le t \le T} E\big(\Vert U_4(t)\Vert^q_{H^{\alpha}_2(\red)}\big) < \infty. 
\eeq 
With (\ref{3.11})-(\ref{3.13}), we finish the proof of the theorem.\hfill\qed 
\medskip 

The next results concern the sample path properties of the solution of (\ref{3.3}). 

\begin{teor} 
\label{t4} 
Let $\sigma$, $b$ be Lipschitz functions satisfying (\ref{3.2}). Fix $k \in \,]0,\infty[$, $\alpha\in [0, k[$ and assume that
there 
exists $\eta\in\,]\frac{\alpha}{k},1[$ such that 
\beq\label{rde48a}
\int_{\red} \frac{\mu(d\xi)}{(1 + |\xi|^2)^{k\eta - \alpha}} < \infty. 
\eeq
Suppose also that $v_0\in H_2^{k \delta + \alpha}(\red)$, for some $\delta\in\,]0,1]$ 
and $\tilde v_0\in H_2^{-(k\gam - \alpha)}(\red)$, 
for some $\gam\in [0,1[$. Set $\theta_0 = \inf(\frac{1}{2}, 1-\eta, \delta, 1-\gam)$. 
Then, for any $q\in [2,\infty[$ and $0\le s\le t\le T$, 
\beqn 
%\label{3.14} 
E\big(\Vert u(t) - u(s)\Vert^q_ {H_2^{\alpha}(\red)}\big) \le C (t - s)^{q\theta}, 
\eeqn 
with $\theta\in\,]0,\theta_0[$. Therefore, the sample paths of the $H^{\alpha}_2(\red)$-valued process $(u(t),\
t\in[0,T])$ solution of (\ref{3.3}) are 
almost surely $\theta$-H\"{o}lder continuous for any $\theta\in\,]0,\theta_0[$. 
%\beqn 
%\theta\in(0, \inf(\frac{1}{2},1-\eta, \delta, 1-\gam)). 
%\eeqn 
\end{teor} 

\noindent{\em Proof.} Fix $0\le s\le t\le 1$. As in the proof of Theorem \ref{t3}, let $U_1(t) = \frac{d}{dt}G(t) * v_0$.
Using the formula $\cos x - \cos y = -2\sin \frac{x+y}{2} \sin \frac{x-y}{2}$, we obtain 
\begin{align} 
\Vert U_1(t) - U_1(s)\Vert^2_{H_2^{\alpha}(\red)}& = 
\Vert (1+|\cdot|^2)^{\frac{\alpha}{2}}(\cos (t|\cdot|^k) - \cos (s|\cdot|^k)) \tf v_0
(\cdot)\Vert^2_{L^2(\red)}\nonumber\\ 
& \le 4 \int_{\red} d\xi\, (1+|\xi|^2)^{\alpha} 
\left(\sin \frac{(t-s) |\xi|^k}{2}\right)^{2\delta} |\tf v_0 (\xi)|^2\nonumber\\ 
& \le (t-s)^{2\delta} \Vert v_0\Vert^2_ {H_2^{k \delta + \alpha}(\red)}\nonumber\\ 
& \le C (t-s)^{2 \delta}. \label{3.14a} 
\end{align} 

Consider now the term $U_2(t) = G(t) * \tilde v_0$. Applying the formula $\sin x - \sin y = 2\cos \frac{x+y}{2} 
\sin\frac{x-y}{2}$ yields 
\begin{align*} 
   \Vert U_2(t) - U_2(s)\Vert^2_{H_2^{\alpha}(\red)} &\leq 
\Vert (1+|\cdot|^2)^{\frac{\alpha}{2}} 
\left(\frac{\sin \frac{(t-s)|\cdot|^k}{2}}{|\cdot|^k}\right)\tf \tilde v_0 (\cdot)\Vert^2_{L^2(\red)}\\ 
& \le T_1 + T_2, 
\end{align*} 
where 
\begin{align*} 
T_1& = (t-s)^2\int_{\red \cap \{|\xi|< 1\}} d \xi\, (1+|\xi|^2)^{\alpha} \frac{|\tf \tilde v_0 (\xi)|^2}{|\xi|^{2k}},\\ 
T_2&= \int_{\red \cap \{|\xi|\ge 1\}} d \xi\, (1+|\xi|^2)^{\alpha} |\tf \tilde v_0 (\xi)|^2 
\frac{\big(\sin \frac{(t-s)|\xi|^k}{2}\big)^{2(1-\gam)}}{|\xi|^{2k}}. 
\end{align*} 
Therefore, 
\begin{align*} 
T_1+T_2& \le C_1 (t-s)^2 \int_{\red} \frac{d\xi}{(1+|\xi|^2)^{k-\alpha}}|\tf \tilde v_0 (\xi)|^2\\ 
& \qquad \qquad + C_2 (t-s)^{2(1-\gam)}\int_{\red} \frac{d\xi}{(1+|\xi|^2)^{k\gam-\alpha}}|\tf \tilde v_0 (\xi)|^2\\ 
& \le C (t-s)^{2(1-\gam)}\Vert \tilde v_0 \Vert^2_{H_2^{-(k\gam-\alpha)}(\red)}. 
\end{align*} 
Consequently, by the assumption on $\tilde v_0$,  
\beq 
\label{3.15} 
\Vert U_2(t) - U_2(s)\Vert_{H_2^{\alpha}(\red)} \le C(t-s)^{1-\gam}. 
\eeq 
Set 
\beqn 
U_3(t) = \int_0^t \int_{\red} G(t-r, \cdot - y) \sig(u(r,y)) M(dr,dy). 
\eeqn 
For any $q\in[2,\infty[$, the following estimate holds: 
\beq 
\label{3.16} 
E\big ( \Vert U_3(t) - U_3(s)\Vert^{q}_ {H^{\alpha}_2(\red)}\big) \le C (t-s) ^{q(\frac{1}{2}\wedge(1-\eta))}. 
\eeq 
Indeed, set $Z(s,y) = \sig(u(s,y))$. Properties (\ref{3.2}) and (\ref{3.9}) imply that 
\beqn 
\sup_{0\le s \le T} E\big(\Vert \sig(u(s))\Vert^{q}_{L^2(\red)}\big) < \infty. 
\eeqn 
Hence, (\ref{3.16}) follows from the upper bound estimate (\ref{2.5}). 

Finally, set 
\beqn 
U_4(t) = \int_0^t ds \int_{\red} dy\, G(t-s, \cdot - y ) b(u(s,y)). 
\eeqn 
Clearly, 
\beqn 
E\big(\Vert U_4(t) - U_4(s) \Vert^{q}_{H_2^{\alpha}(\red)}\big) \le A(s,t) + B(s,t), 
\eeqn 
where 
\begin{align*} 
A(s,t)& = E\big(\Vert \int_s^t dr \int_{\red} dy\, G(t-r, \cdot - y) b(u(r,y))\Vert^{q}_ {H_2^{\alpha}(\red)}\big)\\ 
B(s,t)& = E\big(\Vert \int_0^s dr \int_ {\red} dy\, \big(G(t-r,\cdot - y) - G(s-r, \cdot - y)\big) b(u(r,y))\Vert^{q}_ 
{H_2^{\alpha}(\red)}\big). 
\end{align*} 

The Cauchy-Schwarz inequality, H\"older's inequality, Plancherel's identity, (\ref{3.2}) and the fact that $\alpha < k$ yield 
\begin{align*} 
A(s,t)& = E\Big(\Big( \int_{\red} dx\, \big|\int_s^t dr\,\int_{\red}dy\, (I - \Delta)^{\frac{\alpha}{2}} 
G(t-r,x-y) b(u(r,y))\big|^2\Big)^{\frac{q}{2}}\Big) \\ 
& \le (t-s)^{\frac{q}{2}} E\Big(\Big( \int_s^t dr\,\int_{\red} dx\, \big|\int_{\red} dy\, (I - \Delta)^{\frac{\alpha}{2}} 
G(t-r,x-y) b(u(r,y))\big|^2\Big)^{\frac{q}{2}}\Big)\\ 
&\le (t-s)^{q-1} \int_s^t dr\, E\Big(\Big(\int_{\red} dx\, \big|\int_{\red} dy\, (I - \Delta)^{\frac{\alpha}{2}} 
G(t-r,x-y) b(u(r,y))\big|^2\Big)^{\frac{q}{2}}\Big)\\ 
& = (t-s)^{q-1}\int_s^t dr\, E\Big(\Big(\int_{\red} d\xi\, (1+|\xi|^2)^{\alpha} |\tf G(t-r)(\xi)|^2 |\tf b(u(r))(\xi)|^2 
\Big)^{\frac{q}{2}}\Big)\\ 
&\le C (t-s)^q. 
\end{align*} 
Analogously, by the formula $\sin x - \sin y = 2\cos\frac{x+y}{2}\sin\frac{x-y}{2}$,
\begin{align*} 
B(s,t)& \le C \int_0^s dr E\Big(\big( \int_{\red} dx\, \big|\int_{\red} dy\,
(I-\Delta)^{\frac{\alpha}{2}} 
\big(G(t-r,x-y) - G(s-r,x-y)\big)\\ 
& \qquad \times b(u(r,y))\big|^2\Big)^{\frac{q}{2}}\Big)\\ 
&\le C \int_0^s dr\, E\Big(\big(\int_{\red} d\xi\,
(1+|\xi|^2)^{\alpha}\Big|\frac{\sin\frac{(t-s)|\xi|^k}{2}}{|\xi|^k}\Big|^2  |\tf
b(u(r))(\xi)|^2\big)^{\frac{q}{2}}\Big)\\
 &\le C \int_0^s dr\, E\Big(\big(\int_{\red} d\xi\,
(1+|\xi|^2)^{\alpha}\frac{(\sin\frac{(t-s)|\xi|^k}{2})^{2(1-\eta)}}  {|\xi|^{2k}}\, |\tf
b(u(r))(\xi)|^2\big)^{\frac{q}{2}}\Big)\\
 &\le (t-s)^{q(1-\eta)}\int_0^s dr\, E\Big(\big(
\int_{\red} d\xi\, (1+|\xi|^2)^{\alpha-k\eta}\, |\tf  b(u(r))(\xi)|^2\big)^{\frac{q}{2}}\Big)\\ 
&\le C (t-s)^{q(1-\eta)},
\end{align*} 
because $\alpha - k \eta < 0$. Consequently, 
\beq 
\label{3.17} 
E\big(\Vert U_4(t) - U_4(s) \Vert^{q}_{H_2^{\alpha}(\red)}\big) \le C (t-s)^{q(1-\eta)}. 
\eeq 
The result then follows from (\ref{3.14a})-(\ref{3.17}). \hfill \qed 
\medskip 

We finish this section with a refinement of the previous theorem in the particular case 
of a covariance measure $\Gam$ given by a Riesz kernel. 

\begin{teor} 
\label{t5} 
Fix $k \in\,]0,\infty[$ and $\alpha\in [0,k[$. Let $\sigma$, $b$, $v_0$, $\tilde v_0$, $\delta$ and $\gamma$ be as in
Theorem \ref{t4}. We assume 
that $\Gam(dx) = |x|^{-\beta}$, with $\beta\in\,]0,2(k-\alpha)[$. Set $\theta_1\in\,]0, 
\inf(1-\frac{\beta+2\alpha}{2k}, \delta, 1-\gam)[$. Then, for any $q\in [2,\infty[$, 
$0\le s\le t\le T$, 
\beq 
\label{3.18} 
E\big(\Vert u(t) - u(s)\Vert^q_ {H_2^{\alpha}(\red)}\big) \le C (t - s)^{q\theta}, 
\eeq 
with $\theta\in\,]0,\theta_1[$. 
Therefore, the sample paths of the $H^{\alpha}_2(\red)$-valued process $(u(t),\ t\in[0,T])$ solution of (\ref{3.3}) are
almost surely $\theta$-H\"{o}lder continuous for any 
\beqn 
\theta\in\,\left]0, \inf\left(1-\frac{\beta+2\alpha}{2k}, \delta, 1-\gam\right)\right[. 
\eeqn 
\end{teor} 

\noindent{\it Proof.} We shall use the same notations as in the proof of Theorem \ref{t4}. 
By Theorem \ref{t2.2} with $Z(s,y) = \sigma(u(s,y))$ (see (\ref{2.9})), 
\beq 
\label{3.19} 
E\big(\Vert U_3(t) - U_3(s)\Vert^q_{H_2^{\alpha}(\red)}\big) \le C (t_2-t_1)^{q(1-\frac{\beta+2\alpha}{2k})}. 
\eeq 
It is easy to check that for $\mu(d\xi) = |\xi|^{-d+\beta}$, the condition (\ref{rde48a}) holds in fact for any
$\eta\in\,](2\alpha+\beta)/(2k),1[$. Consequently, (\ref{3.17}) yields 
\beq 
\label{3.20} 
E\big(\Vert U_4(t) - U_4(s)\Vert^q_{H_2^{\alpha}(\red)}\big) \le C (t_2-t_1)^{q\theta_2}, 
\eeq 
for any $\theta_2\in\,]0, 1-(2\alpha+\beta)/(2k)[$. 

The upper bound estimate (\ref{3.18}) is a consequence of (\ref{3.14a}), (\ref{3.15}), 
(\ref{3.19}), (\ref{3.20}) and H\"older continuity of the $H_2^{\alpha}(\red)$-valued process $(u(t),\ t\in[0,T])$ 
follows from Kolmogorov's continuity condition. 
\hfill\qed 
\bigskip

\noindent{\bf Acknowledgement.} The second named author thanks the
{\it Institut de Math\'ematiques} of the {\it Ecole Polytechnique F\'ed\'erale de Lausanne} for 
its hopitality and financial support during a visit where part of this work was carried out.


\begin{thebibliography}{99} 
\bibitem{adler} R. J. Adler: {\it An Introduction to Continuity, Extrema, and Related Topics 
for General Gaussian Processes}. Institute of Mathematical Statistics Lecture Notes-Monographs Series, 
Vol 12, 1990. 
\bibitem{carmonanualart} R. Carmona, D. Nualart: {\it Random nonlinear wave equations: smoothness of the solutions}. 
Probab. Theory Related Fields 79, 469-508 (1988). 
\bibitem{dalangfrangos} R.C. Dalang, N.E. Frangos: {\it The stochastic wave equation in two spatial dimensions}. 
Annals of Probab. 26, 1, 187-212, 1998. 
\bibitem{dalang} R.C. Dalang: {\it Extending the martingale measure stochastic integral with applications 
to spatially homogeneous spde's}. Electronic J. of Probability, Vol 4, 1999. 
\bibitem{dalangmueller} R.C. Dalang, C. Mueller: {\it Some non-linear SPDE's that are second order in 
time}. Electronic J. of Probability, Vol 8, 1, 1-21, 2003. 
\bibitem{dalangss} R.C. Dalang, M. Sanz-Sol\'e: {\it H\"older-Sobolev properties of the solution of the 
stochastic wave equation in dimension three}. In preparation. 
\bibitem{donoghue} W.F. Donoghue: {\it Distributions and Fourier transforms}. Academic Press, New York, 1969. 
\bibitem{gasquetw}C. Gasquet, P. Witomski: {\it Fourier Analysis and Applications}. Texts in Applied 
Mathematics 30. Springer Verlag, 1999. 
\bibitem{hormander} L. H\"{o}rmander: {\it Lectures on Nonlinear Hyperbolic Differential Equations}. Springer Verlag,
1997. 
\bibitem{KZ} A. Karkzewska, J. Zabczyk: ``Stochastic PDE's
with function-valued solutions'', in Cl\'ement Ph., den 
Hollander F., van Neerven J. and de Pagter B. (Eds), ``Infinite-dimensional
stochastic analysis'', Proceedings of the Colloquium of the Royal Netherlands
Academy of Arts and Sciences, 1999, Amsterdam.
\bibitem{krylov} N.V. Krylov: {\it An analytic approach to spde's}. In: {\it Stochastic Partial Differential 
Equations: Six Perspectives} (R.A. Carmona, B. Rozovskii, Eds.), pp. 185-242. Mathematical Surveys and Monographs, 
Vol 64, American Mathematical Society, 1999. 
\bibitem{metivier} M. M\'{e}tivier: {\it Semimartingales, a Course on Stochastic Processes}. de Gruyter Studies in 
Mathematics 2. Walter de Gruyter, 1982. 
\bibitem{milletss} A. Millet, M. Sanz-Sol\'{e}: {\it A stochastic wave equation in two space dimensions: smoothness 
of the law}. Annals of Probab. 27, 803-844, 1999. 
\bibitem{pz} S. Peszat, J. Zabczyk: {\it Nonlinear stochastic wave and heat equations}. Probab. Theory Related Fields,
116, 421-443, 2000.
\bibitem{RY} D. Revuz, M. Yor: {\em Continuous martingales and Brownian motion.} Third edition. Springer-Verlag, Berlin,
1999.
\bibitem{sssarra} M. Sanz-Sol\'{e}, M. Sarr\`a: {\it H\"older continuity for the stochastic heat
equation with spatially correlated noise}. In:  {\it Stochastic analysis, random fields and applications} (R.C. Dalang, M. Dozzi
,
F. Russo Eds.), pp. 259-268, Progress in Probability 52,  Birkh\"auser, Basel, 2002.
\bibitem{sanzsole} M. Sanz-Sol\'{e}: {\it A course on Malliavin Calculus with Applications to Stochastic 
Partial Differential Equations}. Lecture Notes Vol. 2. Institut de Matem\`{a}tica. Universitat de Barcelona, 2004. 
\bibitem{schwartz} L. Schwartz: {\it Th\'eorie des distributions}. Hermann, Paris, 1966. 
\bibitem{shimakura}N. Shimakura: {\it Partial differential operators of elliptic type}. Translations of Mathematical
Monographs, 99. American Mathematical Society, 1992.
\bibitem{sogge} C. D. Sogge: {\it Lectures on Nonlinear Wave Equations}. Monographs in Analysis, Vol II. 
International Press, 1995. 
\bibitem{treves} F. Treves: {\em Basic Linear Partial Differential Equations.} Academic Press, 1975.
\bibitem{walsh}J.B. Walsh: {\it An introduction to stochastic partial 
differential equations}, \'Ecole d'\'et\'e de Probabilit\'es de Saint Flour 
XIV, Lecture Notes in Mathematics, Vol. 1180, Springer Verlag, 1986. 

\end{thebibliography}
\end{document}